%% file: Poisson2.tex
\documentclass[10pt,twoside,a4paper]{article}
\usepackage{makeidx}
\usepackage[dvips]{graphicx}
\usepackage{amscd}
\usepackage{amsmath}
\usepackage{amssymb}
\usepackage{latexsym}
\usepackage[francais]{babel}
\input{xy}
\xyoption{all}

\newcounter{num}[section]
\newcounter{nums}[subsection]
\renewcommand{\thenum}{\arabic{num}}
\renewcommand{\thenums}{\arabic{nums}}
\newcommand{\parags}{\addtocounter{nums}{1}{\thesubsection.\thenums.~}}
\newcommand{\parag}{\addtocounter{num}{1}{\thesection.\thenum.~}}

\renewcommand{\sl}{{\mathfrak{sl}(2)}}
\newcommand{\g}{{\mathfrak{g}}}
\newcommand{\h}{{\mathfrak{h}}}
\newcommand{\C}{{\mathbb{C}}}
\renewcommand{\P}{{\cal P}}
\newcommand{\nequiv}{\equiv \hspace{-3.5mm} /\:}

\newcommand{\sdirect}{\neq \hspace{-2mm} \neq}

\title{Le groupe des traces de Poisson
de la vari\'et\'e quotient $\h\oplus \h^*/W$ en rang $2$}
\author{Jacques Alev \footnote{e-mail : jacques.alev@univ-reims.fr}\:
et Lo{\"\i}c Foissy \footnote{e-mail : loic.foissy@univ-reims.fr}\\
\\
{\small{\it Laboratoire de Math\'ematiques - UMR6056, Universit\'e de Reims}}\\
\small{{\it Moulin de la Housse - BP 1039 - 51687 REIMS Cedex 2, France}}
}
\date{}

\newtheorem{prop}{\indent Proposition}

\newtheorem{lemme}[prop]{\indent Lemme}
\newtheorem{theo}[prop]{\indent Th\'eor\`eme}

\begin{document}

\maketitle 

\input{chap1.tex}

\input{chap2.tex}
\input{chap3.tex}
\input{chap4.tex}

\input{chap5.tex}
\input{chap6.tex}

\input{chap7.tex}
\input{chap8.tex}

\input{chap9.tex}

\bibliographystyle{amsalpha-fr}
\bibliography{biblio}

\end{document}

%% file: chap1.tex
\section{Introduction}

\parag Soit $V$ un espace vectoriel symplectique sur
$\C$, $dim_\C\, V=2l$, et soit $G$ un sous-groupe fini de $Sp(V)$.
Posons $S=\C[V]$.
L'alg\`ebre des fonctions r\'eguli\`eres invariantes $S^G$ 
h\'erite naturellement d'une structure d'alg\`ebre de
Poisson induite et munit ainsi la vari\'et\'e quotient ${\cal
X}=V/G$ d'une structure de vari\'et\'e de Poisson. On peut alors
consid\'erer la d\'eformation non commutative de $\cal X$
d\'efinie par l'alg\`ebre des invariants $A_l(\C)^G$, o\`u
$A_l(\C)$ d\'esigne l'alg\`ebre de Weyl de rang $l$. Il
existe une litt\'erature r\'ecente abondante sur l'\'etude des
d\'esingularisations (symplectiques) de $\cal X$, le calcul des
(co)homologies \'equivariantes de $\cal X$ et les alg\`ebres de
r\'eflexions symplectiques qui en fournissent les d\'eformations
les plus riches (voir \cite{AF,EG,FU,BG}).\\

\parag Il existe deux familles particuli\`eres
d'exemples de la situation d\'ecrite en $1.1$. La premi\`ere
consiste \`a commencer avec un sous-groupe fini $\Gamma$ de
$SL(2,\C)$, \`a consid\'erer $V=(\C^2)^n$, $n\in \mathbb{N}^*$,
et \`a prendre pour $G$ le produit en couronne de $\Gamma$
par $S_n$, produit semi-direct de $\Gamma^n$ par le groupe
sym\'etrique $S_n$. Comme cas particulier, nous avons ici les
surfaces dites de Klein, ${\cal X}_{\Gamma}=\C^2/\Gamma$,
pour $\Gamma$ de type $A_n, D_n, E_{6,7,8}$. La deuxi\`eme famille
consiste \`a commencer avec une alg\`ebre de Lie simple $\g$, une
sous-alg\`ebre de Cartan $\h$, consid\'erer $V=\h\oplus \h^*$ et
prendre pour $G$ le groupe de Weyl $W$ avec l'action diagonale.\\

\parag Conform\'ement \`a l'esprit des d\'eformations
alg\'ebriques, la question standard consiste \`a comparer la
(co)homologie de Poisson de $\cal X$ \`a la (co)homologie de
Hochschild de $A_l(\C)^G$. Le th\'eor\`eme 6.1 de \cite{AFLS}
donne le calcul complet de la (co)homologie de  Hochschild de
$A_l(\C)^G$. Si $a_k$ d\'esigne le nombre de classes de conjugaison
de $G$ agissant dans la repr\'esentation $V$ avec un sous-espace
de points fixes de dimension $k$, $0\leq k\leq 2l$, alors,
$dim_\C\,HH_k(A_l(\C)^G)=a_k$. Curieusement, le calcul
de la (co)homologie de Poisson de $\cal X$ se r\'ev\`ele bien plus
compliqu\'e : cela est d\^u au fait qu'en ce qui concerne
$A_l(\C)^G$, on dispose d'une \'equivalence de Morita qui ram\`ene
les calculs \`a ceux relatifs au produit crois\'e $A_l(\C))\sdirect G$ 
qui se pr\^ete beaucoup mieux aux calculs
(co)homologiques. Pour le calcul de la (co)homologie de Poisson,
nous n'avons actuellement que la m\'ethode directe.\\

\parag Dans la g\'en\'eralit\'e du paragraphe 1.1, on d\'emontre
dans \cite{BEG} que $HP_0({\cal X})$ est de dimension finie.
Le but de cette note est de pr\'esenter le
calcul du groupe d'homologie de Poisson de $\cal X$ en degr\'e
z\'ero, $HP_0({\cal X})$, dans diff\'erents cas o\`u $HP_0({\cal X})$ a
m\^eme dimension que $HH_0(A_l(\C)^G)$.
Cette co{\"\i}ncidence de dimensions peut \^etre 
interpr\'et\'ee comme \'etant le reflet d'une bonne 
d\'eformation, comme c'est le cas pour les 
surfaces de Klein (\cite{AL}).
Pour les trois exemples en rang $2$ de
la deuxi\`eme famille, on trouve ainsi :\\

{\bf Th\'eor\`eme.} Avec les notations pr\'ec\'edentes, on a 
l'\'egalit\'e :
$$dim_\C\,HP_0({\h\oplus \h^*}/W)=dim_\C\,HH_0\left(A_l(\C)^W\right)$$ 
et cette dimension commune vaut $1$ en
type $ A_2$, $2$ en type $B_2$ et $3$ en type $G_2$.
\\

La m\'ethode est la suivante : dans les diff\'erents cas \'etudi\'es,
la composante homog\`ene de degr\'e $2$ de $\C[V]^G$,
munie du crochet de Poisson, est une alg\`ebre de Lie que nous noterons $\g$, 
agissant sur $\C[V]^G$ de mani\`ere semi-simple. De plus, les composantes isotypiques
non triviales de $\C[V]^G$ sous l'action de $\g$ sont inclus dans $\{\g,\C[V]^G\}$ ;
par suite, le calcul de $HP_0({\cal X})$ se restreint \`a calculer la composante isotypique
triviale de $\C[V]^G$ et son intersection avec $\{\C[V]^G,\C[V]^G\}$. 
Dans les exemples de la deuxi\`eme partie (groupes cycliques)
et de la derni\`ere partie (sous-groupes de $(\mathbb{Z}/3\mathbb{Z})^n$),
l'alg\`ebre de Lie $\g$ est ab\'elienne. Dans les exemples de la troisi\`eme partie
(groupes de Weyl de rang 2), il s'agit de $\sl$. 
Dans les exemples de la cinqui\`eme partie (sous-groupes de $(\mathbb{Z}/2\mathbb{Z})^n$),
il s'agit de $\sl^{\oplus n}$.\\

\parag Le papier s'organise de la mani\`ere suivante : nous
consid\'erons d'abord  une famille d'exemples simples, puis nous continuons en
effectuant une \'etude commune des trois groupes de Weyl de rang $2$ en
pr\'esentant les d\'etails d'un calcul bas\'e essentiellement sur
le pl\'ethysme des repr\'esentations de $\sl$ ; 
nous poursuivons par une remarque-question sur le fait que l'id\'eal
d\'eriv\'e de Poisson de $\C[{\cal X}]$, alg\`ebre de
fonctions r\'eguli\`eres sur la vari\'et\'e quotient affine $\cal
X$, est \'egalement un id\'eal associatif. A notre connaissance,
les premiers exemples o\`u cela ne se produit pas sont les cas
$B_2$ et $G_2$. Nous donnons ensuite une pr\'esentation de ces trois alg\`ebres
d'invariants.
La partie suivante expose une famille d'exemples pour lesquels
ces deux dimensions diff\`erent.
Nous terminons par l'\'etude d'une famille d'exemples qui montrent que la diff\'erence
des deux dimensions consid\'er\'ees dans le th\'eor\`eme principal peut \^etre
arbitrairement grande.\\

{\bf Remerciements.} Le premier auteur tient \`a remercier 
Y. Berest, D. Farkas, B. Fu et T. Lambre pour des
conversations fructueuses \`a l'origine de ce travail.

%% file: chap2.tex
\section{Une famille d'exemples simples}

\parag Soit $n \geq 2$ et soit $G=C_n=(\sigma)$ le groupe cyclique d'ordre $n$,
agissant sur $\C^2$ par le caract\`ere $\zeta=e^{2i\pi/n}$ et donc sur
$V=\C^2\oplus (\C^2)^*$. Soit $S=S(V)$ l'alg\`ebre sym\'etrique de $V$,
$A$ l'alg\`ebre de Weyl $A_2(\C)$, avec pour coordonn\'ees respectives
$x_1$, $x_2$, $y_1$, $y_2$ et $p_1$, $p_2$, $q_1$, $q_2$.
L'action  de $G$ sur $S$ et $A$ est alors donn\'ee par :
$$\begin{array}{c|c|c|c|c}
&x_1&x_2&y_1&y_2\\
\hline
\sigma&\zeta x_1&\zeta x_2& \zeta^{n-1} y_1&\zeta^{n-1} y_2
\end{array} \hspace{1cm} \mbox{ et } \hspace{1cm}
\begin{array}{c|c|c|c|c}
&p_1&p_2&q_1&q_2\\
\hline
\sigma&\zeta p_1&\zeta p_2& \zeta^{n-1} q_1&\zeta^{n-1} q_2
\end{array}$$
On note $S^G$ et $A^G$ les alg\`ebres d'invariants respectives. 
Par la structure symplectique standard, $S$ est une alg\`ebre de Poisson
et $G$ agit sur $S$ par automorphismes de Poisson, ce qui implique que $S^G$ 
est une sous-alg\`ebre de Poisson de $S$. Autrement dit, $S^G$ est une alg\`ebre de Poisson
commutative et peut \^etre consid\'er\'ee comme l'alg\`ebre des fonctions r\'eguli\`eres
sur une vari\'et\'e de Poisson alg\'ebrique affine.\\

\parag D'autre part, $S$ est (en tant qu'alg\`ebre de Poisson) gradu\'ee par le degr\'e total,
le crochet de Poisson \'etant homog\`ene de degr\'e $-2$. Comme $G$ agit de mani\`ere homog\`ene
sur $S$, $S^G$ est une sous-alg\`ebre de Poisson gradu\'ee de $S$. 
On note $S^G(n)$ sa composante homog\`ene de degr\'e $n$.
En particulier, $S^G(2)$ munie du crochet de Poisson est une alg\`ebre de Lie
et l'identit\'e de Jacobi implique que $S^G$ est un $S^G(2)$-module.
Remarquons que $x_1y_1$ et $x_2y_2$ appartiennent \`a $S^G(2)$ et que 
$\g=Vect(x_1y_1,x_2y_2)$ est une sous-alg\`ebre de Lie ab\'elienne de $S^G(2)$.
Par suite, $S^G$ est un $\g$-module. 
Pour tout $(\alpha_1,\alpha_2,\beta_1,\beta_2)\in \mathbb{N}^4$ :
\begin{eqnarray*}
\{x_1y_1,x_1^{\alpha_1}y_1^{\beta_1}x_2^{\alpha_2}y_2^{\beta_2}\}
&=&(\beta_1-\alpha_1)x_1^{\alpha_1}y_1^{\beta_1}x_2^{\alpha_2}y_2^{\beta_2},\\
\{x_2y_2,x_1^{\alpha_1}y_1^{\beta_1}x_2^{\alpha_2}y_2^{\beta_2}\}
&=&(\beta_2-\alpha_2)x_1^{\alpha_1}y_1^{\beta_1}x_2^{\alpha_2}y_2^{\beta_2},
\end{eqnarray*}
donc $S$ est une somme directe de $\g$-modules de dimension $1$: $\g$ agit 
de mani\`ere semi-simple sur $S$. D\'ecomposons $S$ en composantes isotypiques :
$$S=\bigoplus_{(i,j)\in \mathbb{Z}^2} S_{(i,j)},$$
avec :
\begin{eqnarray*}
S_{(i,j)}&=&\{X \in S\:/\:\{x_1y_1,X\}=iX,\:\{x_2y_2,X\}=jX\}\\
&=&Vect(x_1^{\alpha_1}y_1^{\beta_1}x_2^{\alpha_2}y_2^{\beta_2}\:/\:
\beta_1-\alpha_1=i,\:\beta_2-\alpha_2=j).
\end{eqnarray*}
De plus, pour tous $(i,j)$, $(k,l) \in \mathbb{Z}^2$, $S_{(i,j)}S_{(k,l)} \subseteq S_{(i+k,j+l)}$ 
et $\{S_{(i,j)},S_{(k,l)}\} \subseteq S_{(i+k,j+l)}$ :
$S$ est ainsi $\mathbb{Z}$-gradu\'ee en tant qu'alg\`ebre de Poisson.\\

\parag D\'ecrivons maintenant $S^G$ :

\begin{prop}
\label{propo1}
\begin{description}
\item[\it i)] Pour tout $(i,j) \in \mathbb{Z}^2$, posons $S^G_{(i,j)}=S^G\cap S_{(i,j)}$.
Alors 
$\displaystyle S^G=\bigoplus_{(i,j)\in \mathbb{Z}^2} S^G_{(i,j)}$.
\item[\it ii)] Si $i+j \nequiv 0[n]$, alors $S^G_{(i,j)}=0$.
\item[\it iii)] Si $i+j \equiv 0[n]$, alors $S^G_{(i,j)}=S_{(i,j)}$.
\item[\it iv)] $S^G_{(0,0)}=\C[t_1,t_2]$, avec $t_1=x_1y_1$ et $t_2=x_2y_2$.
\item[\it v)] $S^G$ est engendr\'ee par $t_1=x_1y_1$, $t_2=x_2y_2$, $x_1^n$, $y_1^n$,
$x_2^n$, $y_2^n$, $x_1^{\alpha} x_2^{n-\alpha}$ $(1\leq \alpha \leq n-1)$,
$y_1^{\alpha} y_2^{n-\alpha}$ $(1\leq \alpha \leq n-1)$, $x_1y_2$ et $x_2y_1$.
\end{description}
\end{prop}

{\bf Preuve.} 
$i)$ Comme $\g \subseteq S^G$, $S^G$ est un sous-$\g$-module de $S$ et donc
se d\'ecompose selon les composantes isotypiques de $S$, d'o\`u le premier point.

$ii)$ et $iii)$ Pour tout $(\alpha_1,\alpha_2,\beta_1,\beta_2)\in \mathbb{N}^4$,
$\sigma.x_1^{\alpha_1}y_1^{\beta_1}x_2^{\alpha_2}y_2^{\beta_2}
=\zeta^{\alpha_1-\beta_1+\alpha_2-\beta_2} x_1^{\alpha_1}y_1^{\beta_1}x_2^{\alpha_2}y_2^{\beta_2}$.
Par suite :
$$S^G=Vect(x_1^{\alpha_1}y_1^{\beta_1}x_2^{\alpha_2}y_2^{\beta_2}\:/\: 
\alpha_1-\beta_1+\alpha_2-\beta_2\equiv 0[n]).$$
Les points $ii)$ et $iii)$ s'en d\'eduisent imm\'ediatement.

$iv)$ On a :
$$S^G_{(0,0)}=Vect\left(x_1^{\alpha_1}y_1^{\beta_1}x_2^{\alpha_2}y_2^{\beta_2}\:/\: 
\beta_1-\alpha_1=\beta_2-\alpha_2=0\right)
=Vect\left(t_1^{\alpha_1}t_2^{\alpha_2}\:/\:(\alpha_1,\alpha_2)\in \mathbb{N}^2\right)
=\C[t_1,t_2].$$

$v)$ Notons $S'$ la sous-alg\`ebre de $S$
engendr\'ee par les \'el\'ements d\'ecrits dans l'\'enonc\'e de la proposition.
De mani\`ere imm\'ediate, ces g\'en\'erateurs propos\'es sont dans $S^G$, donc $S'\subseteq S^G$.

Soit $X=x_1^{\alpha_1}y_1^{\beta_1}x_2^{\alpha_2}y_2^{\beta_2}$, avec  
$\alpha_1-\beta_1+\alpha_2-\beta_2\equiv 0[n]$.
Suivant les valeurs de $\alpha_1,\beta_1$ et de $\alpha_2,\beta_2$, 
on peut alors \'ecrire $X$ sous l'une des formes suivantes :
%\begin{eqnarray*}
%X&=&t_1^{\gamma_1} t_2^{\gamma_2} x_1^{\alpha'_1}x_2^{\alpha'_2}\\
%\mbox{ou }X&=&t_1^{\gamma_1} t_2^{\gamma_2} x_1^{\alpha'_1}y_2^{\beta'_2}\\
%\mbox{ou }X&=&t_1^{\gamma_1} t_2^{\gamma_2} y_1^{\beta'_1}x_2^{\alpha'_2}\\
%\mbox{ou }X&=&t_1^{\gamma_1} t_2^{\gamma_2} y_1^{\beta'_1}y_2^{\beta'_2}.
%\end{eqnarray*}
$$X=t_1^{\gamma_1} t_2^{\gamma_2} x_1^{\alpha'_1}x_2^{\alpha'_2}
\mbox{ ou }X=t_1^{\gamma_1} t_2^{\gamma_2} x_1^{\alpha'_1}y_2^{\beta'_2}
\mbox{ ou }X=t_1^{\gamma_1} t_2^{\gamma_2} y_1^{\beta'_1}x_2^{\alpha'_2}
\mbox{ ou }X=t_1^{\gamma_1} t_2^{\gamma_2} y_1^{\beta'_1}y_2^{\beta'_2}.$$
En effectuant une division euclidienne par $n$, 
on peut \'ecrire $X$ sous l'une des formes suivantes :
\begin{eqnarray*}
X&=&t_1^{\gamma_1} t_2^{\gamma_2} (x_1^n)^{\delta_1} (x_2^n)^{\delta_2} 
x_1^{\alpha''_1}x_2^{\alpha''_2}, \: 0\leq \alpha''_1,\alpha''_2<n\\
\mbox{ou }X&=&t_1^{\gamma_1} t_2^{\gamma_2} (x_1^n)^{\delta_1} (y_2^n)^{\delta_2}
x_1^{\alpha''_1}y_2^{\beta''_2}, \: 0\leq \alpha''_1,\beta''_2<n\\
\mbox{ou }X&=&t_1^{\gamma_1} t_2^{\gamma_2}(y_1^n)^{\delta_1} (x_2^n)^{\delta_2} 
y_1^{\beta''_1}x_2^{\alpha''_2}, \: 0\leq \beta''_1,\alpha''_2<n\\
\mbox{ou }X&=&t_1^{\gamma_1} t_2^{\gamma_2} (y_1^n)^{\delta_1} (y_2^n)^{\delta_2}
y_1^{\beta''_1}y_2^{\beta''_2}, \: 0\leq \beta''_1,\beta''_2<n.
\end{eqnarray*}
Dans le premier cas, on a $\alpha''_1+\alpha''_2\equiv 0[n]$ et donc 
$\alpha''_1+\alpha''_2=0$ ou $n$. Par suite, soit $X$ est de la forme 
$t_1^{\gamma_1} t_2^{\gamma_2} (x_1^n)^{\delta_1} (x_2^n)^{\delta_2}$,
soit $X$ est de la forme $t_1^{\gamma_1} t_2^{\gamma_2} (x_1^n)^{\delta_1} (x_2^n)^{\delta_2} 
\left(x_1^{\alpha''_1}x_2^{n-\alpha''_1}\right)$, avec $0<\alpha''_1<n$. Donc $X \in S'$.
De m\^eme, dans le quatri\`eme cas, Soit $X$  est de la forme 
$t_1^{\gamma_1} t_2^{\gamma_2} (y_1^n)^{\delta_1} (y_2^n)^{\delta_2}$,
soit $X$ est de la forme $t_1^{\gamma_1} t_2^{\gamma_2} (y_1^n)^{\delta_1} (y_2^n)^{\delta_2} 
\left(y_1^{\beta''_1}y_2^{n-\beta''_1}\right)$, avec $0<\beta''_1<n$,
donc $X \in S'$.

Dans le deuxi\`eme cas, on a $\alpha''_1-\beta''_2\equiv 0[n]$, donc $\alpha''_1=\beta''_2$
et $X$ s'\'ecrit :
$$X=t_1^{\gamma_1} t_2^{\gamma_2} (x_1^n)^{\delta_1} (y_2^n)^{\delta_2}(x_1y_2)^{\alpha''_1}.$$
Par suite, $X \in S'$. De m\^eme, dans le troisi\`eme cas, on montre que $X \in S'$ et donc
$S^G=S'$. $\Box$\\

\parag Soit $M$ un sous-$\g$-module simple non trivial de $S^G$. Alors $\g.M$ est un sous-module
non nul de $M$, donc est \'egal \`a $M$. Par suite, $M=\g.M = \{g, M\}\subseteq \{S^G,S^G\}$.
Donc les composantes isotypiques non triviales de $S^G$ sont incluses dans $\{S^G,S^G\}$. 
On en d\'eduit :
\begin{eqnarray*}
\bigoplus_{(i,j) \in \mathbb{Z}^2-\{(0,0)\}} S^G_{(i,j)}&\subseteq & \{S^G,S^G\},\\
S^G&=&S^G_{(0,0)} + \{S^G,S^G\},\\
HP_0(S^G)&=&\frac{S^G_{(0,0)}}{S^G_{(0,0)} \cap \{S^G,S^G\}}.
\end{eqnarray*}
Nous pouvons maintenant montrer que $HP_0(S^G)$
est de dimension  $n-1$.

\begin{theo} 
\label{theo2}
On a l'égalité $dim_{\C}\:HP_0(S^G)=n-1$.
\end{theo}

{\bf Preuve.} Utilisons l'identit\'e g\'en\'erale suivante,
valable pour toute alg\`ebre de Poisson, et qui est 
une application directe de l'identit\'e de Leibniz :
$$\{ab,c\}=\{a,bc\}+\{b,ca\}.$$

En utilisant cette identit\'e autant que n\'ecessaire et le point $v)$
de la proposition \ref{propo1} :
\begin{eqnarray*}
\{S^G,S^G\}&=&
\{t_1, S^G\}+\{t_2, S^G\}+
\{x_1y_2, S^G\}+\{x_2y_1, S^G\}\\
&&+\{x_1^n, S^G\}+\{y_1^n, S^G\}+\{x_2^n, S^G\}+\{y_2^n, S^G\}\\
&&+\sum_{\alpha=1}^{n-1}
\left(\{x_1^{\alpha} x_2^{n-\alpha}, S^G\}
+\{y_1^{\alpha} y_2^{n-\alpha}, S^G\}\right).
\end{eqnarray*}
Par homog\'en\'eit\'e du crochet de Poisson :
\begin{eqnarray*}
\{S^G,S^G\}\cap S_{(0,0)}^G&=&
\{t_1, S^G_{(0,0)}\}+\{t_2, S^G_{(0,0)}\}+
\{x_1y_2, S^G_{(1,-1)}\}+\{x_2y_1, S^G_{(-1,1)}\}\\
&&+\{x_1^n, S^G_{(n,0)}\}+\{y_1^n, S^G_{(-n,0)}\}+\{x_2^n, S^G_{(0,n)}\}+\{y_2^n, S^G_{(0,-n)}\}\\
&&+\sum_{\alpha=1}^{n-1}
\left(\{x_1^{\alpha} x_2^{n-\alpha}, S^G_{(\alpha,n-\alpha)}\}
+\{y_1^{\alpha} y_2^{n-\alpha}, S^G_{(-\alpha,\alpha-n)}\}\right)\\
&=&
\{t_1, S^G_{(0,0)}\}+\{t_2, S^G_{(0,0)}\}+
\{x_1y_2, S^G_{(1,-1)}\}+\{x_2y_1, S^G_{(-1,1)}\}\\
&&+\sum_{\alpha=0}^{n}
\left(\{x_1^{\alpha} x_2^{n-\alpha}, S^G_{(\alpha,n-\alpha)}\}
+\{y_1^{\alpha} y_2^{n-\alpha}, S^G_{(-\alpha,\alpha-n)}\}\right).
\end{eqnarray*}
$S_{(0,0)}
=\C[t_1,t_2]$ \'etant la composante isotypique triviale de $S^G$, 
$\{t_1, S^G_{(0,0)}\}=\{t_2, S^G_{(0,0)}\}=(0)$.
D'autre part,
\begin{eqnarray*}
S^G_{(1,-1)}&=&
Vect(x_1^{\alpha_1}y_1^{\beta_1}x_2^{\alpha_2}y_2^{\beta_2}\:/\:\beta_1-\alpha_1=1,\:
\beta_2-\alpha_2=-1)\\
&=&Vect(t_1^{\alpha_1}t_2^{\beta_2}y_1x_2\:/\: \alpha_1,\beta_2 \in \mathbb{N})\\
&=&y_1x_2\C[t_1,t_2].
\end{eqnarray*}
Calculons :
$$\{x_1y_2, y_1x_2 t_1^{\alpha_1}t_2^{\alpha_2}\}
=(1+\alpha_1)t_1^{\alpha_1}t_2^{\alpha_2+1} -(1+\alpha_2) t_1^{\alpha_1+1}t_2^{\alpha_2}.$$
Par suite :
$$\{x_1y_2, S^G_{(1,-1)}\}
=Vect\left((1+\alpha_1)t_1^{\alpha_1}t_2^{\alpha_2+1} -(1+\alpha_2) t_1^{\alpha_1+1}t_2^{\alpha_2}
\:/\alpha_1,\alpha_2 \in \mathbb{N}\right).$$
Un calcul semblable montre que 
$\{x_2y_1, S^G_{(-1,1)}\}=\{x_1y_2, S^G_{(1,-1)}\}$.

Fixons $0\leq \alpha \leq n$. 
\begin{eqnarray*}
S^G_{(\alpha,n-\alpha)}&=&
Vect\left(x_1^{\alpha_1}y_1^{\beta_1}x_2^{\alpha_2}y_2^{\beta_2}\:/\:
\beta_1-\alpha_1=\alpha,\: \beta_2-\alpha_2=n-\alpha\right)\\
&=&Vect\left(t_1^{\alpha_1}t_2^{\alpha_2}y_1^{\alpha}y_2^{n-\alpha}
\:/\: \alpha_1,\alpha_2 \in \mathbb{N}\right)\\
&=&y_1^{\alpha}y_2^{n-\alpha}\C[t_1,t_2].
\end{eqnarray*}
Calculons :
$$\{x_1^{\alpha} x_2^{n-\alpha},
y_1^{\alpha}y_2^{n-\alpha} t_1^{\alpha_1}t_2^{\alpha_2}\}
=\alpha(\alpha+\alpha_1) t_1^{\alpha-1+\alpha_1}t_2^{n-\alpha+\alpha_2}
+(n-\alpha)(n-\alpha+\alpha_2)t_1^{\alpha+\alpha_1}t_2^{n-\alpha-1+\alpha_2}.$$
En particulier, pour $\alpha=0$ et $\alpha=n$ :
%\begin{eqnarray*}
%\{x_2^n,y_2^nt_1^{\alpha_1}t_2^{\alpha_2}\}
%&=&n(n+\alpha_2)t_1^{\alpha_1}t_2^{n-1+\alpha_2},\\
%\{x_1^n,y_1^n t_1^{\alpha_1}t_2^{\alpha_2}\}
%&=&n(n+\alpha_1) t_1^{n-1+\alpha_1}t_2^{\alpha_2}.
%\end{eqnarray*}
$$\{x_2^n,y_2^nt_1^{\alpha_1}t_2^{\alpha_2}\}
=n(n+\alpha_2)t_1^{\alpha_1}t_2^{n-1+\alpha_2},\hspace{1cm}
\{x_1^n,y_1^n t_1^{\alpha_1}t_2^{\alpha_2}\}
=n(n+\alpha_1) t_1^{n-1+\alpha_1}t_2^{\alpha_2}.$$
Par suite, 
\begin{eqnarray*}
&&\{x_1y_2, S^G_{(1,-1)}\}+\{x_2y_1, S^G_{(-1,1)}\}+\sum_{\alpha=0}^{n}
\{x_1^{\alpha} x_2^{n-\alpha}, S^G_{(\alpha,n-\alpha)}\}\\
&=&Vect\left((1+\alpha_1)t_1^{\alpha_1}t_2^{\alpha_2+1} -(1+\alpha_2) t_1^{\alpha_1+1}t_2^{\alpha_2}
\:/\alpha_1,\alpha_2 \in \mathbb{N}\right)
+t_1^{n-1}\C[t_1,t_2]+t_2^{n-1}\C[t_1,t_2].
\end{eqnarray*}
On obtient un r\'esultat semblable pour $\{y_1^{\alpha} y_2^{n-\alpha}, S^G_{(-\alpha,\alpha-n)}\}$.

On a donc :
\begin{eqnarray*}
&&\{S^G,S^G\}\cap S_{(0,0)}^G\\
&=&Vect\left((1+\alpha_1)t_1^{\alpha_1}t_2^{\alpha_2+1} -(1+\alpha_2) t_1^{\alpha_1+1}t_2^{\alpha_2}
\:/\alpha_1,\alpha_2 \in \mathbb{N}\right)
+t_1^{n-1}\C[t_1,t_2]+t_2^{n-1}\C[t_1,t_2]\\
&=&Vect\left((1+\alpha_1)t_1^{\alpha_1}t_2^{\alpha_2+1} -(1+\alpha_2) t_1^{\alpha_1+1}t_2^{\alpha_2}
\:/\alpha_1+\alpha_2<n-1\right)
+Vect\left(t_1^{\alpha_1}t_2^{\alpha_2} \:/\: \alpha_1+\alpha_2\geq n-1\right).
\end{eqnarray*}
Une base de $S_{(0,0)}^G/ \{S^G,S^G\}\cap S_{(0,0)}^G$
est donc donn\'ee par $\left(\overline{t_1^0},\ldots,\overline{t_1^{n-2}}\right)$,
d'o\`u $S_{(0,0)}^G/ \{S^G,S^G\}\cap S_{(0,0)}^G$ est de dimension $n-1$. $\Box$\\

\parag {\bf Remarques.} 
\begin{description}
\item[\it i)]
 \begin{description}
\item[\it a)] Si $n=2$, $\{S^G,S^G\}$ est l'id\'eal d'augmentation de $S^G$.
\item[\hspace{-4mm}\it b)] Si $n=3$, $S_{(0,0)}^G/ \{S^G,S^G\}\cap S_{(0,0)}^G$ est 
$Vect(t_1-t_2)+Vect\left(t_1^{\alpha_1}t_2^{\alpha_2} \:/\: \alpha_1+\alpha_2\geq n-1\right)$
et est donc l'id\'eal de $S_{(0,0)}^G$ engendr\'e par $t_1-t_2$, $t_1^2$, $t_2^2$ et $t_1t_2$.
\item[\hspace{-4mm}\it c)] Si $n\geq 4$, $t_1-t_2 \in \{S^G,S^G\}\cap S_{(0,0)}^G$
et $(t_1-t_2)t_2 \notin \{S^G,S^G\}\cap S_{(0,0)}^G$, 
donc $\{S^G,S^G\}$ n'est pas un id\'eal de $S^G$.
\end{description}
\item[\it ii)] Si $n \geq 3$, alors $S^G(2)=\g$. Si $n=2$, alors $S^G(2)=S(2)$ est isomorphe
\`a $\sl \oplus \sl$. On peut alors en d\'eduire une preuve plus rapide  du th\'eor\`eme 
\ref{theo2} en utilisant la composante isotypique triviale de $S^G$ sous l'action de 
$\sl \oplus \sl$ plut\^ot que l'action de $\g$. En effet, cette composante
isotypique $S^G_{\sl \oplus \sl}$ v\'erifie :
\begin{eqnarray*}
S^G_{\sl \oplus \sl}&\subseteq& Ker(\{x_1^2,.\}) \cap Ker(\{x_2^2,.\}) 
\cap Ker(\{y_1^2,.\}) \cap Ker(\{y_2^2,.\}) \\
&\subseteq&Ker\left(x_1 \frac{\partial}{\partial y_1}\right)
\cap Ker\left(x_2 \frac{\partial}{\partial y_2}\right)
\cap Ker\left(y_1 \frac{\partial}{\partial x_1}\right)
\cap Ker\left(y_2 \frac{\partial}{\partial x_2}\right)\\
&\subseteq & \C[x_1,x_2,y_2] \cap  \C[x_1,x_2,y_1] \cap  \C[x_2,y_1,y_2] \cap  \C[x_1,y_1,y_2]\\
&\subseteq& \C,
\end{eqnarray*}
donc $S^G_{\sl \oplus \sl}=\C$. Par suite, $S^G=\C + \{S^G,S^G\}$ 
et on  v\'erifie ais\'ement que $\C \cap \{S^G,S^G\}=(0)$.
\end{description}

%% file: chap3.tex
\section{M\'ethode utilis\'ee pour les trois groupes de Weyl de rang $2$}

\label{parag1}

\parag Nous montrerons par la suite que, dans le cas des groupes $A_2$, $B_2$ et $G_2$,
les hypoth\`eses suivantes sont v\'erifi\'ees :\\

{\bf Hypoth\`eses}
\begin{description}
\item[\it a)] $S=S(V)$, avec $V$ de dimension $4$, gradu\'ee avec les \'el\'ements de $V$
homog\`enes de degr\'e $1$. La composante homog\`ene de degr\'e $n$ de $S$ est not\'ee $S(n)$.
De plus, $S$ est munie d'un crochet de Poisson $\{-,-\}$ homog\`ene de degr\'e $-2$.
\item[\it b)] $G$ est un groupe fini agissant par 
automorphismes de  Poisson homog\`enes de degr\'e $0$
sur $S$. On note $S^G$ l'ensemble des \'el\'ements de $S$ invariants sous l'action de $G$ ;
c'est une sous-alg\`ebre de Poisson gradu\'ee de $S$.
\item[\it c)] Il existe trois \'el\'ements non nuls de $S^G(2)$ not\'es $E,F,H$ v\'erifiant :
$$\{E,F\}=H,\hspace{1cm} \{H,E\}=2E, \hspace{1cm} \{H,F\}=-2F.$$
Autrement dit, $Vect(E,F,H)$ muni de $\{-,-\}$ est une alg\`ebre de Lie isomorphe \`a $\sl$.
Alors $\sl$ agit sur $S$ de la mani\`ere suivante : pour tous $P \in S$, $X \in \sl$,
$$X.P=\{X,P\}.$$
Cette action est homog\`ene de degr\'e $0$. Par suite, pour tout $n \in \mathbb{N}$,
$S(n)$ est somme directe d'espaces de poids :
$$S(n)=\bigoplus_{i\in \mathbb{Z}} S(n)_i,$$
o\`u $S(n)_i=\{P \in S(n)\:/\: H.P=iP\}$.
\item[\it d)] $S(1)$ se d\'ecompose en deux sous-espaces de poids $S(1)_1$ et $S(1)_{-1}$,
tous deux de dimension $2$.\\
\end{description}

\parag {\bf Notations}. Soit $M$ un $\sl$-module. 
On note $M_\sl$ sa composante isotypique triviale :
$$M_\sl=\{P \in M\:/\: \forall X \in \sl,\: X.P=0\}.$$
Dans le cas de $S$, comme $E$, $F$ et $H$ sont $G$-invariants, l'action de $\sl$ et l'action
de $G$ commutent et donc :
$$(S^G)_\sl=(S_\sl)^G.$$
Nous noterons par la suite $S_\sl^G$ cette sous-alg\`ebre de Poisson.\\

\parag Soit $M$ un sous-$\sl$-module simple de $S^G$. Si $M$ n'est pas trivial, alors $\sl.M$
est un sous-module non nul de $M$, donc  est \'egal \`a $M$. Par suite, 
$M \subseteq \{S^G,S^G\}$. Les composantes isotypiques non triviales de $S^G$
sont donc incluses dans $\{S^G,S^G\}$ et par suite :
\begin{eqnarray*}
S^G&=&S^G_\sl+\{S^G,S^G\},\\
HP_0(S^G)&=&\displaystyle\frac{S^G_\sl}{\{S^G,S^G \} \cap S^G_\sl}.
\end{eqnarray*}

Nous nous int\'eressons donc maintenant \`a $S^G_\sl$ et \`a $\{S^G,S^G \} \cap S^G_\sl$.

\begin{prop}
\label{propo3}
\begin{description}
\item[\it i)] L'espace vectoriel $S_\sl(2)$ est de dimension $1$. On notera $D$ un g\'en\'erateur fix\'e de cet espace. 
\item[\it ii)] La sous-alg\`ebre $S_\sl$ est engendr\'ee par $D$.
\item[\it iii)] Il existe $N \in \mathbb{N}^*$, tel que $S_\sl^G$ soit la sous-alg\`ebre de $S$
engendr\'ee par $D^N$.
\end{description}
\end{prop}

\begin{lemme}
Les applications $m:S\otimes S \longrightarrow S$ et $\{-,-\}:S\otimes S \longrightarrow S$,
respectivement donn\'ees par le produit et le crochet de Poisson de $S$, 
sont des morphismes de $\sl$-modules.
\end{lemme}

{\bf Preuve.}
Soient $X \in \sl$, $P,Q\in S$.
\begin{eqnarray*}
m(X.(P \otimes Q))&=&m(X.P \otimes Q+P \otimes X.Q)\\
&=&\{X,P\}Q+P\{X,Q\}\\
&=&\{X,PQ\}\\
&=&X.m(P \otimes Q),\\
\\
\{X.(P \otimes Q)\}&=&\{X.P \otimes Q+P \otimes X.Q\}\\
&=&\{\{X,P\},Q\}+\{P,\{X,Q\}\}\\
&=&\{X,\{P,Q\}\}\\
&=&X.\{P \otimes Q\}.
\end{eqnarray*}
(On a utilis\'e l'\'egalit\'e de  Jacobi pour l'avant derni\`ere \'egalit\'e).
Donc $m$ et $\{-,-\}$ sont des morphismes de $\sl$-modules. $\Box$\\

{\bf Preuve de la proposition \ref{propo3}}.
Graduons $S$ sur $\mathbb{N}^2$ en mettant les \'el\'ements de 
$S(1)_1$ homog\`enes de degr\'e $(1,0)$
et les \'el\'ements de $S(1)_{-1}$ homog\`enes de degr\'e $(0,1)$. On note $S(i,j)$
les composantes homog\`enes de $S$ pour cette graduation. 
Alors la s\'erie formelle de Poincar\'e-Hilbert $\Phi(x,y)$ de $S$ est :
$$\Phi(x,y)=\sum_{i,j} dim_\C \: S(i,j)\: x^i y^j=\left(\frac{1}{1-x}\right)^2
\left(\frac{1}{1-y}\right)^2
=\sum_{i,j} (i+1)(j+1)x^iy^j.$$

Comme $m$ est un morphisme de $\sl$-modules, $m$ est homog\`ene pour le poids et le degr\'e. 
Autrement dit, pour tous $i,j \in \mathbb{Z}$, $m,n,j \in \mathbb{N}$, 
$S(m)_i S(n)_j \subseteq S(m+n)_{i+j}$. D'autre part, 
on remarque que $S(1,0)=S(1)_1$ et $S(0,1)=S(1)_{-1}$.
Comme $S$ est engendr\'ee par 
$S(1)=S(1,0)\oplus S(0,1)$, si $i,j \in \mathbb{N}$ :
$$S(i,j)=S(1,0)^i S(0,1)^j =S(1)_1^i S(1)_{-1}^j \subseteq S(i+j)_{i-j}.$$
On en d\'eduit :
\begin{equation}
\label{E1}
S(n)_k=\displaystyle\bigoplus_{i+j=n,i-j=k}S(i,j).
\end{equation}
On note $\chi_n$ le caract\`ere du $\sl$-module $S(n)$ :
$$\chi_n(q)=\sum_{k\in \mathbb{Z}} dim_\C \:S(n)_k\: q^k.$$
On pose alors :
$$\chi(q,h)=\sum_{n=0}^{\infty} \chi_n(q)h^n
=\sum_{n,k} dim_\C \:S(n)_k\:h^n q^k.$$
Par (\ref{E1}) :
\begin{eqnarray*}
\chi(q,h)&=&\sum_{k\in \mathbb{Z},n\in \mathbb{N}} dim_\C \: S(n)_k h^nq^k\\
&=& \sum_{k\in \mathbb{Z},n\in \mathbb{N}}\:\sum_{i+j=n,\:i-j=k} dim_\C \: S(i,j) h^nq^k\\
&=&\sum_{i,j \in \mathbb{N}} dim_\C \: S(i,j) h^{i+j}q^{i-j}\\
&=&\sum_{i,j \in \mathbb{N}} dim_\C \: S(i,j) (hq)^i\left(\frac{h}{q}\right)^j\\
&=&\Phi(hq,h/q)\\
&=&\sum_{i,j} (i+1)(j+1)h^{i+j}q^{i-j}.
\end{eqnarray*}
Comme la dimension de $S(n)_\sl$ est la diff\'erence 
du terme constant et du terme en $q^2$ de $\chi_n$,
cherchons ces deux coefficients.
Il faut donc prendre :
\begin{enumerate} 
\item Pour le terme constant de $\chi_n$ :
$$ \left\{\begin{array}{rcl}
  i+j&=&n \\
  \\[-1mm]
  i-j&=&0,
\end{array}\right.
\Longleftrightarrow
  \left\{\begin{array}{rcl}
  i&=&\displaystyle \frac{n}{2} \\
  \\[-1mm]
  j&=&\displaystyle \frac{n}{2}.
\end{array}\right.$$
\item Pour le terme en $q^2$ de $\chi_n$ :
$$\left\{\begin{array}{rcl}
  i+j&=&n \\
  \\[-1mm]
  i-j&=&2,
\end{array}\right.
\Longleftrightarrow
\left\{\begin{array}{rcl}
  i&=&\displaystyle \frac{n}{2}+1 \\
  \\[-1mm]
  j&=&\displaystyle \frac{n}{2}-1,
\end{array}\right.
 $$
\end{enumerate}
Par suite, si $n$ est impair, ces deux coefficients sont nuls 
et donc $S(n)_\sl=(0)$. Si $n=2l$ est pair, le terme constant est $(l+1)^2$
et le terme en $q^2$ est $l(l+2)$, donc :
$$dim_\C\:S(n)_\sl=(l+1)^2-l(l+2)=1.$$
En particulier, $S(2)_\sl$ est de dimension $1$, ce qui implique le premier point.
Comme $m$ est un morphisme de $\sl$-modules, pour tout $l\in \mathbb{N}$,
$D^l$ est un \'el\'ement non nul de $S(2l)_\sl$ et donc forme une base de $S(2l)_\sl$.
Par suite, $S_\sl=\C[D]$. 

Comme $G$ agit de mani\`ere homog\`ene sur $S_\sl$
et que $S_\sl(2)$ est de dimension $1$, il existe un unique caract\`ere $\kappa$ de $G$ 
tel que pour tout $\sigma \in G$,
$\sigma.D=\kappa(\sigma)D$.
Comme $S_\sl^G$ est une sous-alg\`ebre gradu\'ee de $S_\sl$ :
\begin{eqnarray*}
S_\sl^G&=&Vect(D^k\:/\:\forall \sigma \in G,\: \sigma.D^k=D^k)\\
&=&Vect(D^k\:/\: \forall \sigma \in G,\: \kappa(\sigma)^k=1)\\
&=&Vect(D^k\:/\: \kappa^k=1)\\
&=&Vect(D^k\:/\: o(\kappa)\mid k).
\end{eqnarray*}
(Comme $G$ est fini, son groupe de caract\`eres est fini et donc $o(\kappa)$ est fini).
Posons $N=o(\kappa)$. Alors $S_\sl=Vect(D^k\:/\: N\mid k)=K[D^N]$. $\Box$\\

\parag Par homog\'en\'eit\'e, 
les composantes homog\`enes de $S_\sl$ \'etant nulles ou de dimension $1$,
$$\{S^G,S^G \} \cap S^G_\sl=Vect\left(D^{kN}\:/\:D^{kN}\in \{S^G,S^G \}\right).$$
D\'ecrivons maintenant un proc\'ed\'e permettant de construire des \'el\'ements 
de $\{S^G,S^G \} \cap S^G_\sl$ :

\begin{prop}
\label{propo5}
Soit $P \in S^G$ un vecteur de plus haut poids $\beta$ homog\`ene de degr\'e $\alpha_1$
et $Q \in S^G$ un vecteur de plus haut poids $\beta$ homog\`ene de degr\'e $\alpha_2$.
On pose $\alpha=(\alpha_1+\alpha_2-2)/2$, qu'on supposera entier. 
On d\'efinit par r\'ecurrence sur $i$ :
$$P^{(\beta)}=P,\hspace{1cm} P^{(\beta-2i)}=\{F,P^{(\beta-2i+2)}\},$$
$$Q^{(\beta)}=Q,\hspace{1cm} Q^{(\beta-2i)}=\{F,Q^{(\beta-2i+2)}\}.$$
Alors il existe deux scalaires $\lambda,\mu\in \C$ tels que pour tout $k \in \mathbb{N}$,
$$\sum_{i=0}^\beta (-1)^i \{P^{(\beta-2i)},Q^{-(\beta-2i)}D^k\}=(\lambda+k\mu)D^{\alpha+k}.$$
En particulier, si $N \mid \alpha+k$,
$(\lambda+k\mu)D^{\alpha+k} \in \{S^G\otimes S^G\}\cap S^G_\sl$.
\end{prop}

{\bf Preuve.}

 Montrons d'abord que l'\'el\'ement suivant est dans $(S\otimes S)_\sl$ :
$$\P=\sum_{i=0}^\beta (-1)^i P^{(\beta-2i)} \otimes Q^{-(\beta-2i)}.$$
Comme $P$ et $Q$ sont des vecteurs de plus haut poids $\beta$, pour tout $i$,
$P^{(\beta-2i)}$ et $Q^{(\beta-2i)}$ sont des vecteurs de poids $\beta-2i$
et $\{F,P^{(-\beta)}\}=\{F,Q^{(-\beta)}\}=0$.
On a donc :
\begin{eqnarray*}
H.\P&=&\sum_{i=0}^\beta (-1)^i \{H,P^{(\beta-2i)}\} \otimes Q^{-(\beta-2i)}
+\sum_{i=0}^\beta (-1)^i P^{(\beta-2i)} \otimes \{H,Q^{-(\beta-2i)}\}\\
&=&\sum_{i=0}^\beta (-1)^i (\beta-2i)P^{(\beta-2i)} \otimes Q^{-(\beta-2i)}
-\sum_{i=0}^\beta (-1)^i (\beta-2i)P^{(\beta-2i)} \otimes Q^{-(\beta-2i)}\\
&=&0,\\
\\
F.\P&=&\sum_{i=0}^\beta (-1)^i \{F,P^{(\beta-2i)}\} \otimes Q^{-(\beta-2i)}
+\sum_{i=0}^\beta (-1)^i P^{(\beta-2i)} \otimes \{F,Q^{-(\beta-2i)}\}\\
&=&\sum_{i=0}^{\beta-1} (-1)^i P^{(\beta-2i-2)} \otimes Q^{-(\beta-2i)}
+\sum_{i=1}^\beta (-1)^i P^{(\beta-2i)} \otimes Q^{-(\beta-2i+2)}\\
&=&0,
\end{eqnarray*}
donc $\P$ est un vecteur de plus bas poids $0$, donc est un \'el\'ement de $(S \otimes S)_\sl$.
Par suite, pour tout $k \in \mathbb{N}$, $\P \otimes D^k \in (S \otimes S \otimes S)_\sl$.

Comme $m$ et $\{-,-\}$ sont des morphismes de $\sl$-modules,
$\{-,-\} \circ (Id \otimes m)(\P \otimes D^k) \in S_\sl$. Par homog\'en\'eit\'e,
il s'agit d'un \'el\'ement homog\`ene de degr\'e $\alpha_1+\alpha_2+2k-2=2(\alpha+k)$, 
donc il existe $\lambda_k\in \C$, $\{-,-\} \circ (Id \otimes m)(\P \otimes D^k)
=\lambda_k D^{\alpha+k}$.
En particulier, posons $\lambda=\lambda_0$.

De la m\^eme mani\`ere, $D \otimes \P \in (S \otimes S \otimes S)_\sl$, homog\`ene de degr\'e
$2\alpha+4$, donc en appliquant $m\circ (\{-,-\} \otimes Id)$, il existe $\mu \in \C$,
$$\sum_{i=0}^\beta (-1)^i \{D,P^{(\beta-2i)}\}Q^{-(\beta-2i)}=-\mu D^{\alpha+1}.$$
D'autre part,
\begin{eqnarray*}
\{-,-\} \circ (Id \otimes m)(\P \otimes D^k)&=&
\sum_{i=0}^\beta (-1)^i \{P^{(\beta-2i)},Q^{-(\beta-2i)}D^k\}\\
&=&\sum_{i=0}^\beta (-1)^i \{P^{(\beta-2i)},Q^{-(\beta-2i)}\}D^k\\
&&+\sum_{i=0}^\beta (-1)^i \{P^{(\beta-2i)},D\}kQ^{-(\beta-2i)}D^{k-1}\\
&=&\lambda D^\alpha D^k+k\mu D^{\alpha+1}D^{k-1}\\
&=&(\lambda+k\mu)D^{\alpha+k},
\end{eqnarray*}
donc $\lambda_k=\lambda+k\mu$. $\Box$\\

\parag {\bf Remarques.}
\begin{description}
\item[\it i)] Les scalaires $\lambda$ et $\mu$ sont d\'etermin\'es par :
\begin{eqnarray*}
\sum_{i=0}^\beta (-1)^i \{P^{(\beta-2i)},Q^{-(\beta-2i)}\}&=&\lambda D^\alpha,\\
\sum_{i=0}^\beta (-1)^i \{P^{(\beta-2i)},D\}Q^{-(\beta-2i)}&=&\mu D^{\alpha+1}.
\end{eqnarray*}
\item[\it ii)] On peut \'eventuellement prendre $P=Q$.
\item[\it iii)] Pour $k=0$, l'\'el\'ement $\sum(-1)^i \{P^{(\beta-2i)},Q^{-(\beta-2i)}\}$ 
peut \^etre repr\'esent\'e \`a l'aide d'arbres de la mani\`ere suivante, 
en utilisant l'anti-sym\'etrie de $\{-,-\}$ :
$$
\begin{picture}(105,70)(-50,0)
\put(0,0){\line(0,0){10}}
\put(0,10){\line(1,1){20}}
\put(22,32){.}
\put(24,34){.}
\put(26,36){.}
\put(30,40){\line(1,1){10}}
\put(40,50){\line(1,1){10}}
\put(-50,60){\line(1,-1){50}}
\put(-30,60){\line(1,-1){40}}
\put(-10,60){\line(1,-1){30}}
\put(10,60){\line(1,-1){20}}
\put(30,60){\line(1,-1){10}}
\put(-53,63){\tiny $P$}
\put(-33,63){\tiny $F$}
\put(-13,63){\tiny $F$}
\put(7,63){\tiny $F$}
\put(27,63){\tiny $F$}
\put(47,63){\tiny $Q$}
\end{picture}
+
\begin{picture}(105,70)(-50,0)
\put(0,0){\line(0,0){10}}
\put(0,10){\line(1,1){20}}
\put(22,32){.}
\put(24,34){.}
\put(26,36){.}
\put(30,40){\line(1,1){10}}
\put(40,50){\line(1,1){10}}
\put(-50,60){\line(1,-1){50}}
\put(-30,60){\line(-1,-1){10}}
\put(-10,60){\line(1,-1){30}}
\put(10,60){\line(1,-1){20}}
\put(30,60){\line(1,-1){10}}
\put(-53,63){\tiny $P$}
\put(-33,63){\tiny $F$}
\put(-13,63){\tiny $F$}
\put(7,63){\tiny $F$}
\put(27,63){\tiny $F$}
\put(47,63){\tiny $Q$}
\end{picture}
+
\begin{picture}(105,70)(-50,0)
\put(0,0){\line(0,0){10}}
\put(0,10){\line(1,1){20}}
\put(22,32){.}
\put(24,34){.}
\put(26,36){.}
\put(30,40){\line(1,1){10}}
\put(40,50){\line(1,1){10}}
\put(-50,60){\line(1,-1){50}}
\put(-30,60){\line(-1,-1){10}}
\put(-10,60){\line(-1,-1){20}}
\put(10,60){\line(1,-1){20}}
\put(30,60){\line(1,-1){10}}
\put(-53,63){\tiny $P$}
\put(-33,63){\tiny $F$}
\put(-13,63){\tiny $F$}
\put(7,63){\tiny $F$}
\put(27,63){\tiny $F$}
\put(47,63){\tiny $Q$}
\end{picture}
+\ldots
$$
$$
\ldots+
\begin{picture}(105,80)(-50,0)
\put(0,0){\line(0,0){10}}
\put(0,10){\line(-1,1){20}}
\put(-24,32){.}
\put(-26,34){.}
\put(-28,36){.}
\put(-30,40){\line(-1,1){10}}
\put(-40,50){\line(-1,1){10}}
\put(50,60){\line(-1,-1){50}}
\put(-30,60){\line(-1,-1){10}}
\put(-10,60){\line(-1,-1){20}}
\put(10,60){\line(-1,-1){30}}
\put(30,60){\line(1,-1){10}}
\put(-53,63){\tiny $P$}
\put(-33,63){\tiny $F$}
\put(-13,63){\tiny $F$}
\put(7,63){\tiny $F$}
\put(27,63){\tiny $F$}
\put(47,63){\tiny $Q$}
\end{picture}
+
\begin{picture}(105,80)(-50,0)
\put(0,0){\line(0,0){10}}
\put(0,10){\line(-1,1){20}}
\put(-24,32){.}
\put(-26,34){.}
\put(-28,36){.}
\put(-30,40){\line(-1,1){10}}
\put(-40,50){\line(-1,1){10}}
\put(50,60){\line(-1,-1){50}}
\put(-30,60){\line(-1,-1){10}}
\put(-10,60){\line(-1,-1){20}}
\put(10,60){\line(-1,-1){30}}
\put(30,60){\line(-1,-1){40}}
\put(-53,63){\tiny $P$}
\put(-33,63){\tiny $F$}
\put(-13,63){\tiny $F$}
\put(7,63){\tiny $F$}
\put(27,63){\tiny $F$}
\put(47,63){\tiny $Q$}
\end{picture}
$$
\end{description}

\parag 
Nous avons donc un proc\'ed\'e permettant de montrer que certains $D^k$ appartiennent
\`a $\{S^G , S^G\}$.
Donnons maintenant un crit\`ere permettant de montrer que $D^k \notin  \{S^G , S^G\}$ :
\begin{prop}
\label{propo6}
Soit $k\in \mathbb{N}$. On suppose que :
$$\left\{
\left(\sum_{i=0}^{k+1}S^G(2k+2-i) \otimes S^G(i)\right)_\sl\right\}=0.$$
Alors $D^k \notin \{S^G , S^G\}$.
\end{prop}

{\bf Preuve.} Supposons $D^k \in \{S^G , S^G\}$. Comme $\{-,-\}$
est un morphisme de $\sl$-modules, $\{-,-\}$ envoie les composantes isotypiques de $S^G \otimes S^G$
sur les composantes isotypiques de $S$ correspondantes, donc
$D^k$ poss\`ede un ant\'ec\'edent dans $(S^G \otimes S^G)_\sl$. Par homog\'en\'eit\'e,
il poss\`ede alors un ant\'ec\'edent dans $(S^G \otimes S^G)_\sl(2k+2)$. 
Par antisym\'etrie de $\{-,-\}$, il poss\`ede un ant\'ec\'edent dans :
$$ \sum_{i=0}^{k+1}(S^G(2k+2-i) \otimes S^G(i))_\sl
=\left(\sum_{i=0}^{k+1}S^G(2k+2-i) \otimes S^G(i)\right)_\sl.$$
Donc l'image par $\{-,-\}$ de ce sous-espace de $S^G\otimes S^G$ est non nulle. $\Box$

%% file: chap4.tex
\section{Calculs pour $B_2$}

\label{sectB2} 

\parag Fixons tout d'abord les notations.
Soit $S=\C[x_1,x_2,y_1,y_2]$, muni du crochet de Poisson donn\'e par :
$$\begin{array}{c|c|c|c|c}
\{-,-\}&x_1&y_1&x_2&y_2\\
\hline
x_1&0&1&0&0\\
\hline
y_1&-1&0&0&0\\
\hline
x_2&0&0&0&1\\
\hline
y_2&0&0&-1&0
\end{array}$$

Le groupe $G=B_2=(\pm 1)^2 \rtimes S_2$ agit par automorphismes de Poisson homog\`enes 
sur l'alg\`ebre $S$ de la mani\`ere suivante :
$$\begin{array}{c|c|c|c|c}
&x_1&y_1&x_2&y_2\\
\hline
(\epsilon,\epsilon',id)&\epsilon x_1&\epsilon y_1&\epsilon' x_2&\epsilon' y_2\\
\hline
(\epsilon,\epsilon',(12))&\epsilon' x_2&\epsilon' y_2&\epsilon x_1&\epsilon y_1\\
\end{array}$$

\begin{lemme}
Consid\'erons les \'el\'ements de $S$ suivants :
$$E=\frac{x_1^2+x_2^2}{2},\hspace{1cm}
F=-\frac{y_1^2+y_2^2}{2},\hspace{1cm}
H=-(x_1y_1+x_2y_2).$$
Alors $(E,F,H)$ v\'erifie l'hypoth\`ese $c)$.
\end{lemme}
{\bf Preuve.} Calculs directs. $\Box$\\

Les \'el\'ements $E$, $F$ et $H$ agissent par d\'erivation sur $S$. Donnons leur action sur les g\'en\'erateurs :
$$\begin{array}{c|c|c|c|c}
&x_1&y_1&x_2&y_2\\
\hline
E&0&x_1&0&x_2\\
\hline
F&y_1&0&y_2&0\\
\hline
H&x_1&-y_1&x_2&-y_2\\
\end{array}$$
Par suite, l'action de $E$, $F$ et $H$ est donn\'ee par :
\begin{eqnarray*}
\{E,-\}&=&x_1\frac{\partial}{\partial y_1}+x_2\frac{\partial}{\partial y_2},\\
\{F,-\}&=&y_1\frac{\partial}{\partial x_1}+y_2\frac{\partial}{\partial x_2},\\
\{H,-\}&=&x_1\frac{\partial}{\partial x_1}+x_2\frac{\partial}{\partial x_2}
-y_1\frac{\partial}{\partial y_1}-y_2\frac{\partial}{\partial y_2}.
\end{eqnarray*}

On a donc $S(1)_{1}=Vect(x_1,x_2)$ et $S(1)_{-1}=Vect(y_1,y_2)$,
donc l'hypoth\`ese $d)$ est satisfaite.\\

\parag On consid\`ere l'\'el\'ement suivant de $S$ :
$$D=\left|
\begin{array}{cc}
x_1 & x_2\\
y_1 & y_2
\end{array} \right|=x_1y_2-y_1x_2.$$
Alors $\{E,D\}=\{F,D\}=\{H,D\}=0$.
D'apr\`es la proposition \ref{propo3},  $S_\sl=\C[D]$. De plus, 
$$(\epsilon_1,\epsilon_2,\sigma).D=\epsilon_1\epsilon_2\varepsilon(\sigma)D,$$
donc, avec les notations de la proposition \ref{propo3}, $N=2$
et $S^{B_2}_\sl=\C[D^2]$.\\

\parag Nous pouvons maintenant d\'emontrer le th\'eor\`eme suivant :
\begin{theo}
On a l'\'egalit\'e $dim_\C\:HP_0(S^{B_2})=2$.
\end{theo}

{\bf Preuve.} On consid\`ere les \'el\'ements suivants :
\begin{eqnarray*}
P&=&x_1^4+x_2^4,\\
Q&=&x_1^4x_2y_2+x_1y_1x_2^4-x_1^3y_1x_2^2-x_1^2x_2^3y_2.
\end{eqnarray*}
Alors $P,Q \in S^{B_2}$. De plus :
%\begin{eqnarray*}
%\{E,P\}&=&0,\\
%\{H,P\}&=&4P,\\
%\\
%\{E,Q\}&=&0,\\
%\{H,Q\}&=&4Q.
%\end{eqnarray*}
$$\{E,P\}=0,\hspace{1cm}
\{H,P\}=4P,\hspace{1cm}
\{E,Q\}=0,\hspace{1cm}
\{H,Q\}=4Q.$$
Donc $P$ est un vecteur de plus haut poids $4$ homog\`ene de degr\'e $4$
et $Q$ est un vecteur de plus haut poids $4$ homog\`ene de degr\'e $6$.
Appliquons la proposition \ref{propo5}. Par un calcul direct, on trouve :
\begin{eqnarray*}
\sum_{i=0}^4 (-1)^i \{P^{(\beta-2i)},Q^{-(\beta-2i)}\}&=&(-48\times 6) D^4,\\
\sum_{i=0}^4 (-1)^i \{P^{(\beta-2i)},D\}Q^{-(\beta-2i)}&=&(-48) D^5.
\end{eqnarray*}
Donc pour tout $k \in \mathbb{N}$, $-48(6+2k)D^{4+2k} \in \{S^{B_2} , S^{B_2}\}$.
Autrement dit, pour tout $l$ pair sup\'erieur ou \'egal \`a $4$, $D^l \in \{S^{B_2} , S^{B_2}\}$.\\

 Reste \`a \'etudier le cas de $1$ et $D^2$. Comme dans le paragraphe \ref{parag1},
graduons $S$ sur $\mathbb{N}^2$ en mettant $x_1$ et $x_2$ en degr\'e $(1,0)$
et $y_1$ et $y_2$ en degr\'e $(0,1)$. Soient $S(i,j)$ les composantes homog\`enes pour
cette graduation. On pose :
$$\Phi(x,y)=\sum_{i,j} dim_\C \:S^{B_2}(i,j)\: x^iy^j.$$

Posons $H=(\pm 1)^2$ et $K=S_2$. 
Comme ${B_2}=H \rtimes K$, $S^{B_2}=(S^{H})^{K}$.
Une base de $S^{H}$ est donn\'ee par les mon\^omes 
$x_1^{\alpha_1}y_1^{\beta_1}x_2^{\alpha_2}y_2^{\beta_2}$, $\alpha_1,\beta_1$ ayant la m\^eme
parit\'e,  $\alpha_2,\beta_2$ ayant la m\^eme parit\'e.
Soit $\xi_{i,j}$ le caract\`ere du $S_2$-module $S^{H}(i,j)$ et posons :
$$\xi(x,y)=\sum_{i,j} \xi_{i,j} x^i y^j.$$
Comme le groupe $S_2$ agit par permutation sur les mon\^omes, $\xi_{i,j}(\sigma)$
est le nombre de mon\^omes de $S^{H}$ de degr\'e $(i,j)$ fix\'es par $\sigma$.
En cons\'equence :
\begin{description}
\item[\it a)] Pour $\sigma=id$, les mon\^omes fix\'es par $id$ sont les mon\^omes 
$x_1^{\alpha_1}y_1^{\beta_1}x_2^{\alpha_2}y_2^{\beta_2}$, $\alpha_1,\beta_1$ ayant la m\^eme
parit\'e,  $\alpha_2,\beta_2$ ayant la m\^eme parit\'e.
Donc $\xi(id)$ est la s\'erie de Poincar\'e-Hilbert de l'alg\`ebre $\mathbb{N}^2$-gradu\'ee
$\mathbb{C}[x_1^2,y_1^2,x_1y_1,x_2y_2]$ :
$$\xi(id)=\left( \frac{1+xy}{(1-x^2)(1-y^2)}\right)^2.$$
\item[\it b)] Pour $\sigma=(12)$, les mon\^omes fix\'es par $(12)$ sont les mon\^omes 
$(x_1x_2)^{\alpha}(y_1y_2)^{\beta}$, $\alpha,\beta$ ayant la m\^eme parit\'e.
Donc $\xi((12))$ est la s\'erie de Poincar\'e-Hilbert de l'alg\`ebre $\mathbb{N}^2$-gradu\'ee
$\mathbb{C}[(x_1x_2)^2,(y_1y_2)^2,x_1y_1x_2y_2]$ :
$$\xi((12))=\frac{1+x^2y^2}{(1-x^4)(1-y^4)}.$$
\end{description}
Par suite, comme $\displaystyle dim_\C(S^H(i,j))=\frac{\xi_{i,j}(id)+\xi_{i,j}((12))}{2}$ :
$$\Phi(x,y)=\frac{\xi(id)+\xi((12))}{2}
=\frac{1+xy+2x^2y^2+xy^3+x^3y+x^3y^3+x^4y^4}{(1+x^2)(1-x)^2(1+x)^2(1+y^2)(1-y)^2(1+y)^2}.$$
On note $\chi_n$ le caract\`ere du $\sl$-module $S^{B_2}(n)$.
On pose alors :
$$\chi(q,h)=\sum_{n=0}^{\infty} \chi_n(q)h^n
=\sum_{n,k} dim_\C \: S^{B_2}(n)_k \:h^n q^k.$$
De mani\`ere semblable \`a la preuve de la proposition \ref{propo3} :

$$\chi(q,h)=\Phi(hq,h/q)
=\frac{q^6+h^2q^6+h^4q^4+2h^4q^6+h^4q^8+h^6q^6+h^8q^6}
{(h^2q^2+1)(h^2+q^2)(hq+1)^2(q-h)^2(h+q)^2(hq-1)^2}.$$
En d\'eveloppant cette fraction rationnelle en s\'erie relativement \`a $h$ :
\begin{eqnarray*}
\chi_0(q)&=&1,\\
\chi_2(q)&=&q^2+1+q^{-2},\\
\chi_4(q)&=&2q^4+2q^2+3+2q^{-2}+2q^{-4},\\
\chi_6(q)&=&2q^6+3q^4+4q^2+4+4q^{-2}+3q^{-4}+2q^{-6}.
\end{eqnarray*}
(Notons que $S^H(n)=(0)$ si $n$ est impair, donc $S^{B_2}(n)=0$ si $n$ est impair).\\

Appliquons la proposition \ref{propo6} pour $k=0$.
Le caract\`ere de $S^{B_2}(2)\otimes S^{B_2}(0)+S^{B_2}(1)\otimes S^{B_2}(1)$
est $\chi_2(q)$, donc $\left(S^{B_2}(2)\otimes S^{B_2}(0)
+S^{B_2}(1)\otimes S^{B_2}(1)\right)_\sl=(0)$,
donc $1 \notin \{S^{B_2} , S^{B_2}\}$.\\

Appliquons la proposition \ref{propo6} pour $k=2$.
Le caract\`ere de $S^{B_2}(6)\otimes S^{B_2}(0)+S^{B_2}(5) \otimes S^{B_2}(1)
+S^{B_2}(4)\otimes S^{B_2}(2)+S^{B_2}(3) \otimes S^{B_2}(3)$
est :
$$\chi_6(q)+\chi_4(q)\chi_2(q)=
4q^6+7q^4+11q^2+11+11q^{-2}+7q^{-4}+4q^{-6},$$
donc $\left(S^{B_2}(6)\otimes S^{B_2}(0)+S^{B_2}(5) \otimes S^{B_2}(1)
+S^{B_2}(4)\otimes S^{B_2}(2)+S^{B_2}(3) \otimes S^{B_2}(3)\right)_\sl=(0)$,
donc $D^2 \notin \{S^{B_2} , S^{B_2}\}$.\\

En conclusion, seuls les polyn\^omes $1$ et $D^2$ n'appartiennent pas au sous-espace $\{S^{B_2} , S^{B_2}\}$.
Donc $dim_\C\:HP_0(S^{B_2})=2$. $\Box$

%% file: chap5.tex
\section{Calculs pour $A_2$ et $G_2$}

\subsection{Calculs pour $A_2$}

\label{sectA2}

\parags Soit $S'=\C[x_1,x_2,x_3,y_1,y_2,y_3]$, muni du crochet de Poisson donn\'e par :
$$\begin{array}{c|c|c|c|c|c|c}
\{-,-\}&x_1&y_1&x_2&y_2&x_3&y_3\\
\hline
x_1&0&1&0&0&0&0\\
\hline
y_1&-1&0&0&0&0&0\\
\hline
x_2&0&0&0&1&0&0\\
\hline
y_2&0&0&-1&0&0&0\\
\hline
x_3&0&0&0&0&0&1\\
\hline
y_3&0&0&0&0&-1&0
\end{array}$$
Le groupe $G=A_2=S_3$ agit par permutation des indices. 
$S$ est la sous-alg\`ebre de $S'$ engendr\'ee par $x_1-x_2$, $x_1-x_3$,
$y_1-y_2$ et $y_1-y_3$. C'est une sous-alg\`ebre de Poisson gradu\'ee et
un sous ${A_2}$-module ; en fait :
$$S'=\C[x_1+x_2+x_3,y_1+y_2+y_3]\otimes S.$$
Comme $x_1+x_2+x_3$ et $y_1+y_2+y_3$ sont $A_2$-invariants :
$$(S')^{A_2}=\C[x_1+x_2+x_3,y_1+y_2+y_3]\otimes S^{A_2}.$$
De plus, $S$ est l'alg\`ebre des polyn\^omes en les \'el\'ements :
\begin{eqnarray*}
a_1&=&2x_1-x_2-x_3,\\
a_2&=&-x_1+2x_2-x_3, \\
b_1&=&2y_1-y_2-y_3,\\
b_2&=&-y_1+2y_2-y_3.
\end{eqnarray*}
On pose \'egalement :
$$a_3=-x_1-x_2+2x_3,\hspace{1cm} b_3=-y_1-y_2+2y_3,$$
de sorte que :
\begin{description}
\item[\it a)] le groupe ${A_2}$ agit sur $a_1$, $a_2$ et $a_3$ par permutation des indices ;
\item[\it b)] le groupe ${A_2}$ agit sur $b_1$, $b_2$ et $b_3$ par permutation des indices ;
\item[\it c)] on a les \'egalit\'es $a_1+a_2+a_3=0$ et $b_1+b_2+b_3=0$.\\
\end{description}
Le crochet de Poisson de $S$ est donn\'e par le tableau suivant :
$$\begin{array}{c|c|c|c|c}
\{-,-\}&a_1&b_1&a_2&b_2\\
\hline
a_1&0&6&0&-3\\
\hline
b_1&-6&0&3&0\\
\hline
a_2&0&-3&0&6\\
\hline
b_2&3&0&-6&0
\end{array}$$

\parags Mettons en \'evidence la copie de $\sl$ de $S^{A_2}$ :

\begin{lemme}
Consid\'erons les \'el\'ements de $S$ suivants :
\begin{eqnarray*}
E&=&\frac{a_1^2+a_2^2+a_3^2}{18}\:=\:\frac{a_1^2+a_2^2+a_1a_2}{9},\\
F&=&-\frac{b_1^2+b_2^2+b_3^2}{18}\:=\:-\frac{b_1^2+b_2^2+b_1b_2}{9},\\
H&=&-\frac{a_1b_1+a_2b_2+a_3b_3}{9}\:=\:-\frac{2a_1b_1+a_1b_2+b_1a_2+2a_2b_2}{9},\\
\end{eqnarray*}
Alors $(E,F,H)$ v\'erifie l'hypoth\`ese $c)$.
\end{lemme}

{\bf Preuve.} 
Comme ${A_2}$ agit par permutation des indices, $E$, $F$ et $H$ sont ${A_2}$-invariants.
Le reste se d\'emontre par des calculs directs. $\Box$\\

Les \'el\'ements $E$, $F$ et $H$ agissent par d\'erivation sur $S$. Donnons leur action sur les g\'en\'erateurs :
$$\begin{array}{c|c|c|c|c}
&a_1&b_1&a_2&b_2\\
\hline
E&0&a_1&0&a_2\\
\hline
F&b_1&0&b_2&0\\
\hline
H&a_1&-b_1&a_2&-b_2\\
\end{array}$$
Par suite, l'action de $E$, $F$ et $H$ est donn\'ee par :
\begin{eqnarray*}
\{E,-\}&=&a_1\frac{\partial}{\partial b_1}+a_2\frac{\partial}{\partial b_2},\\
\{F,-\}&=&b_1\frac{\partial}{\partial a_1}+b_2\frac{\partial}{\partial a_2},\\
\{H,-\}&=&a_1\frac{\partial}{\partial a_1}+a_2\frac{\partial}{\partial a_2}
-b_1\frac{\partial}{\partial b_1}-b_2\frac{\partial}{\partial b_2}.
\end{eqnarray*}

D'autre part, on a donc $S(1)_{+1}=Vect(a_1,a_2)$ et $S(1)_{-1}=Vect(b_1,b_2)$,
donc l'hypoth\`ese $d)$ est satisfaite.\\

\parags On consid\`ere l'\'el\'ement suivant de $S$ :
$$D=\frac{1}{27}\left|
\begin{array}{ccc}
1&1&1\\
x_1 & x_2 & x_3\\
y_1 & y_2 & y_3
\end{array} \right|=a_1b_2-b_1a_2.$$
Alors $\{E,D\}=\{F,D\}=\{H,D\}=0$.
D'apr\`es la proposition \ref{propo3},  $S_\sl=\C[D]$. De plus, 
$\sigma.D=\varepsilon(\sigma)D$,
donc, avec les notations de la proposition \ref{propo3}, $N=2$.\\

\parags Nous pouvons maintenant d\'emontrer le th\'eor\`eme suivant :
\begin{theo}
On a l'\'egalit\'e 
$dim_\C\:HP_0(S^{A_2})=1$.
\end{theo}

{\bf Preuve.}  On consid\`ere l'\'el\'ement suivant :
$$P=-\frac{a_1^3+a_2^3+a_3^3}{3}=a_1^2a_2+a_1a_2^2.$$
Comme le groupe $A_2$ agit par permutation des indices, $P \in S^{A_2}$. De plus,
$\{E,P\}=0$ et $\{H,P\}=3$, donc $P$ est un vecteur de plus haut poids $3$,
homog\`ene de degr\'e $3$. Utilisons la proposition \ref{propo5}, avec $P=Q$.
Par un calcul direct, on trouve :
\begin{eqnarray*}
\sum_{i=0}^3 (-1)^i \{P^{(\beta-2i)},P^{-(\beta-2i)}\}&=&(-36\times 4) D^2,\\
\sum_{i=0}^3 (-1)^i \{P^{(\beta-2i)},D\}P^{-(\beta-2i)}&=&(-36) D^3.
\end{eqnarray*}
Donc pour tout $k \in \mathbb{N}$, $-36(4+2k)D^{2+2k} \in \{S^{A_2} , S^{A_2}\}$.
Autrement dit, pour tout $l$ pair sup\'erieur ou \'egal \`a $2$, $D^l \in \{S^{A_2} , S^{A_2}\}$.\\

Reste \`a \'etudier le cas de $1$. 
On proc\`ede comme pour $B_2$, avec des notations semblables.
Comme dans le paragraphe \ref{parag1},
$S'$ est gradu\'ee en mettant $x_1,x_2$ en degr\'e $(1,0)$ et $y_1,y_2$ en degr\'e $(0,1)$.
On pose : 
$$\Phi'(x,y)=\sum_{i,j} dim_\C \:(S')^{A_2}(i,j)\: x^iy^j,\:
\Phi(x,y)=\sum_{i,j} dim_\C \:S^{A_2}(i,j)\: x^iy^j.$$

Une base de $S'$ est donn\'ee par les mon\^omes 
$x_1^{\alpha_1}y_1^{\beta_1}x_2^{\alpha_2}y_2^{\beta_2}x_3^{\alpha_3}y_3^{\beta_3}$.
Soit $\xi'_{i,j}$ le caract\`ere du ${A_2}$-module $S'(i,j)$ et posons :
$$\xi'(x,y)=\sum_{i,j} \xi'_{i,j} x^i y^j.$$
Comme ${A_2}$ agit par permutation sur les mon\^omes, $\xi'_{i,j}(\sigma)$
est le nombre de mon\^omes de $S^{A_2}$ de degr\'e $(i,j)$ fix\'es par $\sigma$.
En cons\'equence :
\begin{description}
\item[\it a)] Pour $\sigma=id$, les mon\^omes fix\'es par $id$ sont les mon\^omes 
$x_1^{\alpha_1}y_1^{\beta_1}x_2^{\alpha_2}y_2^{\beta_2}x_3^{\alpha_3}y_3^{\beta_3}$.
Donc $\xi'(id)$ est la s\'erie de Poincar\'e-Hilbert de l'alg\`ebre $\mathbb{N}^2$-gradu\'ee
$\mathbb{C}[x_1,x_2,x_3,y_1,y_2,y_3]$ :
$$\xi'(id)=\left(\frac{1}{(1-x)(1-y)}\right)^3.$$
\item[\it b)] Pour $\sigma=(12)$, les mon\^omes fix\'es par $(12)$ sont les mon\^omes 
$(x_1x_2)^{\alpha}(y_1y_2)^{\beta}x_3^{\alpha_3}y_3^{\beta_3}$. 
Donc $\xi'((12))$ est la s\'erie de Poincar\'e-Hilbert de l'alg\`ebre $\mathbb{N}^2$-gradu\'ee
$\mathbb{C}[x_1x_2,x_3,y_1y_2,y_3]$ :
$$\xi'((12))=\frac{1}{(1-x^2)(1-y^2)(1-x)(1-y)}.$$
\item[\it c)] Pour $\sigma=(123)$, les mon\^omes fix\'es par $(123)$ sont les mon\^omes 
$(x_1x_2x_3)^{\alpha}(y_1y_2y_3)^{\beta}$. 
Donc $\xi'((123))$ est la s\'erie de Poincar\'e-Hilbert de l'alg\`ebre $\mathbb{N}^2$-gradu\'ee
$\mathbb{C}[x_1x_2x_3,y_1y_2y_3]$ :
$$\xi'((123))=\frac{1}{(1-x^3)(1-y^3)}.$$
\end{description}
Par suite, comme 
$\displaystyle dim_\C(S'(i,j))=\frac{\xi'_{i,j}(id)+3\xi'_{i,j}((12))+2\xi'_{i,j}((123))}{6}$, on a :
$$\Phi'(x,y)=\frac{\xi'(id)+3\xi'((12))+2\xi'((123))}{6}.$$
Comme $(S')^{A_2}=S[x_1+x_2+x_3,y_1+y_2+y_3]\otimes S^{A_2}$, 
$\displaystyle \Phi'=\frac{1}{(1-x)(1-y)}\Phi$ et donc :
$$\Phi(x,y)=\frac{1+xy+xy^2+x^2y+x^2y^2+x^3y^3}{(x+1)(x^2+x+1)(1-x)^2(y+1)(y^2+y+1)(1-y)^2}.$$
Comme pour $B_2$, en notant $\chi_n$ le caract\`ere du $\sl$-module $S^{A_2}(n)$, on obtient :
$$\sum_{n=0}^\infty \chi_n(q)h^n=\Phi(hq,h/q)
=\frac{q^5+h^2q^5+h^3q^4+h^3q^6+h^4q^5+h^6q^5}
{(hq+1)(h+q)(h^2+hq+q^2)(h^2q^2+hq+1)(h-q)^2(1-hq)^2}.$$
En d\'eveloppant cette fraction rationnelle en s\'erie relativement \`a $h$, il vient :
\begin{eqnarray*}
\chi_0(q)&=&1,\\
\chi_1(q)&=&0,\\
\chi_2(q)&=&q^2+1+q^{-2},\\
\chi_3(q)&=&q^3+q+q^{-1}+q^{-3},\\
\chi_4(q)&=&q^4+q^2+2+q^{-2}+q^{-4},\\
\chi_5(q)&=&q^5+2q^3+2q+2q^{-1}+2q^{-3}+q^{-5},\\
\chi_6(q)&=&2q^6+2q^4+3q^2+3+3q^{-2}+2q^{-4}+2q^{-6},\\
\chi_7(q)&=&q^7+2q^5+3q^3+3q+3q^{-1}+3q^{-3}+2q^{-5}+q^{-7},\\
\chi_8(q)&=&2q^8+3q^6+4q^4+4q^2+5+4q^{-2}+4q^{-4}+3q^{-6}+2q^{-8},\\
\chi_9(q)&=&2q^9+3q^7+4q^5+5q^3+5q+5q^{-1}+5q^{-3}+4q^{-5}+3q^{-7}+2q^{-9},\\
\chi_{10}(q)&=&2q^{10}+3q^{8}+5q^{6}+5q^{4}+6q^2+6+6q^{-2}+5q^{-4}+5q^{-6}
+3q^{-8}+2q^{-10}.
\end{eqnarray*}
Appliquons la proposition \ref{propo6} pour $k=0$.
Le caract\`ere de $S^{A_2}(2)\otimes S^{A_2}(0)+S^{A_2}(1)\otimes S^{A_2}(1)$
est $\chi_2(q)$, donc $\left(S^{A_2}(2)\otimes S^{A_2}(0)
+S^{A_2}(1)\otimes S^{A_2}(1)\right)_\sl=(0)$,
d'o\`u $1 \notin \{S^{A_2} , S^{A_2}\}$.\\

En conclusion, seul le polyn\^ome $1$  n'appartient pas \`a $\{S^{A_2} , S^{A_2}\}$.
Donc $dim_\C\:HP_0(S^{A_2})=1$. $\Box$

%% file: chap6.tex
\subsection{Calculs pour $G_2$}

\label{sectG2}
\parags On reprend les notations du paragraphe pr\'ec\'edent.
Comme $G_2=A_2\times (\pm Id)$, $S^{G_2}=\left(S^{A_2}\right)^{(\pm Id)}$ et par suite :
$$S^{G_2}=\bigoplus_{k=0}^{\infty} S^{A_2}(2k).$$
Donc $E,F,H,D^2\in S^{G_2}$. Par suite, $S^{G_2}_\sl=\C[D^2]$.\\

\parags Nous pouvons maintenant d\'emontrer le th\'eor\`eme suivant :
\begin{theo}
On a l'\'egalit\'e $dim_\C\:HP_0(S^{G_2})=3$.
\end{theo}

{\bf Preuve.}  On consid\`ere les \'el\'ements suivants de $S$ :
\begin{eqnarray*}
P&=&a_1^6+a_2^6+a_3^6\\
&=&2a_1^6 +2 a_2^6 +6 a_1^5 a_2 +15 a_1^4 a_2^2 +20 a_1^3 a_2^3 
+15 a_1^2a_2^4 +6 a_1a_2^5,\\
Q&=&a_1^6(a_2b_2+a_3b_3)+a_2^6(a_1b_1+a_3b_3)+a_3^6(a_1b_1+a_2b_2)\\
&&-a_1^5b_1(a_2^2+a_3^2)-a_2^5b_2(a_1^2+a_3^2)-a_3^5b_3(a_1^2+a_2^2)\\
&=&-2 a_1^6 a_2 b_2 -2 a_1 b_1 a_2^6 -5 a_1^5 a_2^2 b_2 +2 a_1^5 b_1 a_2^2 +5 a_1^3 a_2^4 b_2 
+5 a_1^4 b_1 a_2^3 +2 a_1^2 a_2^5 b_2 -5 a_1^2 b_1 a_2^5.
\end{eqnarray*}
Comme le groupe $A_2$ agit par permutation sur les indices, $P$ et $Q$ sont $A_2$-invariants.
Comme ils sont homog\`enes de degr\'es pairs, ils sont dans $S^{G_2}$.
De plus :
$$\{E,P\}=0,\hspace{1cm}
\{H,P\}=6P,\hspace{1cm}
\{E,Q\}=0,\hspace{1cm}
\{H,Q\}=6Q,$$
donc $P$ et $Q$ sont des vecteurs de plus haut poids $6$. 
Utilisons la proposition \ref{propo5}.
Par un calcul direct, on trouve :
\begin{eqnarray*}
\sum_{i=0}^6 (-1)^i \{P^{(\beta-2i)},Q^{-(\beta-2i)}\}&=&(25\:920\times 8) D^6,\\
\sum_{i=0}^6 (-1)^i \{P^{(\beta-2i)},D\}Q^{-(\beta-2i)}&=&(25\:920) D^7.
\end{eqnarray*}
Donc pour tout $k \in \mathbb{N}$, $25\:920(8+2k)D^{6+2k} \in \{S^{G_2}, S^{G_2}\}$.
Autrement dit, pour tout $l$ pair sup\'erieur ou \'egal \`a $6$, $D^l \in \{S^{G_2},S^{G_2}\}$.\\

Reste \`a \'etudier les cas de $1$, $D^2$ et $D^4$. 
Comme $1 \notin \{S^{A_2}, S^{A_2}\}$, $1\notin \{S^{G_2}, S^{G_2}\}$.
Appliquons la proposition \ref{propo6} pour $k=4$.
 Nous avons calcul\'e le caract\`ere de $S^{A_2}(n)$
pour $n \leq 10$ dans la section pr\'ec\'edente et donc le caract\`ere 
$\chi_n(q)$ de $S^{G_2}(n)$ :
\begin{eqnarray*}
\chi_0(q)&=&1,\\
\chi_2(q)&=&q^2+1+q^{-2},\\
\chi_4(q)&=&q^4+q^2+2+q^{-2}+q^{-4},\\
\chi_6(q)&=&2q^6+2q^4+3q^2+3+3q^{-2}+2q^{-4}+2q^{-6},\\
\chi_8(q)&=&2q^8+3q^6+4q^4+4q^2+5+4q^{-2}+4q^{-4}+3q^{-6}+2q^{-8},\\
\chi_{10}(q)&=&2q^{10}+3q^{8}+5q^{6}+5q^{4}+6q^2+6+6q^{-2}+5q^{-4}+5q^{-6}
+3q^{-8}+2q^{-10}.
\end{eqnarray*}
(Si $n$ est impair, $\chi_n(q)=0$).
Par suite, le caract\`ere de $S^{G_2}(6) \otimes S^{G_2}(0)+S^{G_2}(4) \otimes S^{G_2}(2)$
est :
$$\chi_6+\chi_2\chi_4=3q^6+4q^4+7q^2+7+7q^{-2}+4q^{-4}+3q^{-6},$$
donc $\left(S^{G_2}(6) \otimes S^{G_2}(0)+S^{G_2}(4) \otimes S^{G_2}(3)\right)_\sl=(0)$. 
Par suite, $D^2 \notin \{S^{G_2}, S^{G_2}\}$.\\

 Appliquons la proposition \ref{propo6} pour $k=4$.
Le caract\`ere de $S^{G_2}(10) \otimes S^{G_2}(0)
+S^{G_2}(8) \otimes S^{G_2}(2)+S^{G_2}(6) \otimes S^{G_2}(4)$ est :
$$\chi_{10}(q)+\chi_8(q)\chi_2(q)+\chi_6(q)\chi_4(q)$$
$$=6q^{10}+12q^8+23q^6+28q^4+35q^2+35+35q^{-2}+28q^{-4}+23q^{-6}+12q^{-8}+6q^{-10},$$
donc $\left(S^{G_2}(10) \otimes S^{G_2}(0)
+S^{G_2}(8) \otimes S^{G_2}(2)+S^{G_2}(6) \otimes S^{G_2}(4) \right)_\sl=(0)$.
Par suite, $D^4 \notin \{S^{G_2}, S^{G_2}\}$.\\

Enfin, seuls les polyn\^omes $1$, $D^2$ et $D^4$ n'appartiennent pas \`a $\{S^{G_2}, S^{G_2}\}$.
Donc $dim_\C\:HP_0(S^{G_2})=3$. $\Box$ \\

\parags {\bf Remarques.}
\begin{description}
\item[\it a)] Le sous-espace $\{S^{A_2}, S^{A_2}\}$ est un id\'eal : c'est l'id\'eal d'augmentation
de $S^{A_2}$.
\item[\it b)] Le sous-espace $\{S^{B_2}, S^{B_2}\}$ n'est pas un id\'eal de $S^{B_2}$. 
En effet, $E=\frac{1}{2}\{H,E\}$, $F=-\frac{1}{2}\{H,F\}$ et $H=\{E,F\}$, donc $E$, $F$ et $H$
appartiennent \`a $\{S^{B_2}, S^{B_2}\}$. Un calcul direct montre que :
$$H^2+4EF=D^2.$$
Or $D^2\notin \{S^{B_2}, S^{B_2}\}$, ce qui montre que ceci n'est pas un id\'eal.
\item[\it c)] Le sous-espace $\{S^{G_2}, S^{G_2}\}$ n'est pas un id\'eal de $S^{G_2}$. En effet,
$E=\frac{1}{2}\{H,E\}$, $F=-\frac{1}{2}\{H,F\}$ et $H=\{E,F\}$, donc $E$, $F$ et $H$
appartiennent \`a $\{S^{G_2}, S^{G_2}\}$. Un calcul direct montre que :
$$H^2+4EF=\frac{-1}{27}D^2.$$
Or $D^2\notin \{S^{G_2}, S^{G_2}\}$, ce qui montre que ceci n'est pas un id\'eal.
\end{description}

%% file: chap7.tex
\section{Pr\'esentations des alg\`ebres d'invariants par g\'en\'erateurs et relations
pour les trois groupes de Weyl de rang 2}

\subsection{Cas de $A_2$}

\parags Reprenons les notations du paragraphe \ref{sectA2}.
Rappelons que $S^{A_2}$ est $\mathbb{N}^2$-gradu\'ee en mettant
$a_1$, $a_2$ homog\`enes de degr\'e $(1,0)$ et
$b_1$, $b_2$ homog\`enes de degr\'e $(0,1)$. La s\'erie formelle de Poincar\'e-Hilbert
de $S^{A_2}$ est :
$$\Phi(x,y)=\frac{1+xy+xy^2+x^2y+x^2y^2+x^3y^3}{(1-x^2)(1-x^3)(1-y^2)(1-y^3)}.$$

\parags Consid\'erons les \'el\'ements suivants de $S^{A_2}$ :
\begin{eqnarray*}
S_1&=&-\frac{a_1a_2+a_2a_3+a_1a_3}{9}=\frac{a_1^2+a_2^2+a_1a_2}{9},\\
T_1&=&-\frac{a_1a_2a_3}{9}=\frac{a_1a_2^2+a_2a_1^2}{9},\\
U_1&=&-\frac{a_1b_1^2+a_2b_2^2+a_3b_3^2}{9}=\frac{2a_1b_1b_2+2a_2b_1b_2+a_1b_2^2+a_2b_1^2}{9},\\
S_2&=&-\frac{b_1b_2+b_2b_3+b_1b_3}{9}=\frac{b_1^2+b_2^2+b_1b_2}{9},\\
T_2&=&-\frac{b_1b_2b_3}{9}=\frac{b_1b_2^2+b_2b_1^2}{9},\\
U_2&=&-\frac{a_1^2b_1+a_2^2b_2+a_3^2b_3}{9}=\frac{2a_1a_2b_1+2a_1a_2b_1+a_1^2b_2+a_2^2b_1}{9},\\
H&=&-\frac{a_1b_1+a_2b_2+a_3b_3}{9}=-\frac{2a_1b_1+a_1b_2+a_2b_1+2a_2b_2}{9}.
\end{eqnarray*}
Comme $A_2=S_3$ agit sur $(a_1,a_2,a_3)$ et $(b_1,b_2,b_3)$ par
permutation des indices, ces \'el\'ements sont bien dans $S^{A_2}$.
On note $R=\C[S_1,T_1,S_2,T_2]$. Par la th\'eorie des fonctions sym\'etriques,
$S_1$ et $T_1$ sont alg\'ebriquement ind\'ependants
et $S_2$ et $T_2$ sont alg\'ebriquement ind\'ependants.
Par suite, $R$ est une alg\`ebre de polyn\^omes \`a quatre variables et sa s\'erie de
Poincar\'e-Hilbert est :
$$\frac{1}{(1-x^2)(1-x^3)(1-y^2)(1-y^3)}.$$
\begin{lemme}
Tout \'el\'ement de $S^{A_2}$ s'\'ecrit de mani\`ere unique :
$$P_0+P_1H+P_2H^2+P_3H^3+Q_1U_1+Q_2U_2,$$
avec $P_0,P_1,P_2,P_3,Q_1,Q_2$ des \'el\'ements de $R$.
Autrement dit, $S^{A_2}$ est un $R$-module libre de base $(1,H,H^2,H^3,U_1,U_2)$.
\end{lemme}

{\bf Preuve.} 
{\it Unicit\'e.} Supposons que $Q\in S^{A_2}$ s'\'ecrive :
$$Q=P_0+P_1H+P_2H^2+P_3H^3+Q_1U_1+Q_2U_2,$$
avec $P_0,P_1,P_2,P_3,Q_1,Q_2$ des \'el\'ements de $R$. On introduit l'action de $A_2$ 
par automorphismes d'alg\`ebre sur $S$ suivante : pour tout $\sigma \in S_3$,
$$\sigma \star a_i=\sigma.a_i=a_{\sigma(i)},\hspace{1cm} \sigma\star b_i=b_i.$$
Comme $S_1$, $T_1$, $S_2$ et $T_2$ sont invariants sous cette action,
les \'el\'ements de $R$ sont invariants sous cette action. Par suite, pour tout $\sigma \in S_3$ :
$$\sigma \star Q=P_0+P_1(\sigma \star H)+P_2(\sigma \star H^2)+P_3(\sigma \star H^3)
+Q_1(\sigma \star U_1)+Q_2(\sigma \star U_2).$$
On a donc :
$$\left(\begin{array}{cccccc}
1&H&H^2&H^3&U_1&U_2\\
(12)\star 1&(12)\star H&(12)\star H^2&(12)\star H^3&(12)\star U_1&(12)\star U_2\\
(23)\star 1&(23)\star H&(23)\star H^2&(23)\star H^3&(23)\star U_1&(23)\star U_2\\
(123)\star 1&(123)\star H&(123)\star H^2&(123)\star H^3&(123)\star U_1&(123)\star U_2\\
(132)\star 1&(132)\star H&(132)\star H^2&(132)\star H^3&(132)\star U_1&(132)\star U_2\\
(13)\star 1&(13)\star H&(13)\star H^2&(13)\star H^3&(13)\star U_1&(13)\star U_2
\end{array}\right)
\left(\begin{array}{c}
P_0\\
P_1\\
P_2\\
P_3\\
Q_1\\
Q_2\\
\end{array}
\right)
=\left(\begin{array}{c}
Q\\
(12)\star Q\\
(123)\star Q\\
(132)\star Q\\
(13)\star Q
\end{array}
\right).$$
Le d\'eterminant de la matrice du membre de gauche est :
\begin{eqnarray*}
&&\frac{-1}{9^7}(a_1-a_2)^3(a_1-a_3)^3(a_2-a_3)^3
(b_1-b_2)^3(b_1-b_3)^3(b_2-b_3)^3\\
&=&\frac{-1}{9^7}(a_1-a_2)^3(2a_1+a_2)^3(a_1+2a_2)^3
(b_1-b_2)^3(2b_1+b_2)^3(b_1+2b_2)^3.
\end{eqnarray*}
Il est donc non nul. Par suite, les $P_i$ et les $Q_j$ 
sont enti\`erement d\'etermin\'es en multipliant
le vecteur du membre de droite par l'inverse de la matrice
du membre de gauche.\\

{\it Existence.} D'apr\`es ce qui pr\'ec\`ede, le $R$-module $M$ engendr\'e par
$1$, $H$, $H^2$, $H^3$, $U_1$ et $U_2$ est libre de base $(1,H,H^2,H^3,U_1,U_2)$.
Par suite, sa s\'erie de Poincar\'e-Hilbert est :
$$\frac{1+xy+xy^2+x^2y+x^2y^2+x^3y^3}{(1-x^2)(1-x^3)(1-y^2)(1-y^3)}.$$
C'est \'egalement la s\'erie de Poincar\'e-Hilbert de $S^{A_2}$.
Comme $M\subseteq S^{A_2}$, $M=S^{A_2}$. $\Box$\\

\parags {\bf Remarque.} La preuve du lemme pr\'ec\'edent est algorithmique
et permet donc de calculer explicitement la d\'ecomposition des \'el\'ements
de $S^{A_2}$ dans la base $(1,H,H^2,H^3,U_1,U_2)$. \\

\parags Nous pouvons maintenant donner une pr\'esentation de $S^{A_2}$ :
\begin{theo}
$S^{A_2}$ est engendr\'ee par $S_1$, $T_1$, $U_1$, $S_2$, $T_2$, $U_2$ et $H$ et les relations
suivantes :
\begin{eqnarray*}
H^4&=&-4S_1^2S_2^2-3T_1T_2H+5S_1S_2H^2-T_1S_2U_1-S_1T_2U_2,\\
HU_1&=&-3S_1T_2-S_2U_2,\\
HU_2&=&-3T_1S_2-S_1U_1,\\
U_1^2&=&12S_1S_2^2-3S_2H^2+3T_2U_2,\\
U_2^2&=&12S_1^2S_2-3S_1H^2+3T_1U_1,\\
U_1U_2&=&9T_1T_2-12S_1S_2H+3H^2.
\end{eqnarray*}
\end{theo}

{\bf Preuve.} D'apr\`es le lemme pr\'ec\'edent, $S^{A_2}$ est engendr\'ee par
$S_1$, $T_1$, $U_1$, $S_2$, $T_2$, $U_2$ et $H$. Les relations entre ces \'el\'ements
sont donn\'ees en d\'ecomposant les produits des \'el\'ements de la base $(1,H,H^2,H^3,U_1,U_2)$
du R-module $S^{A_2}$ dans cette m\^eme base.
Des calculs directs utilisant l'algorithme de la preuve du lemme pr\'ec\'edent 
donnent le r\'esultat annonc\'e. $\Box$\\

\parags Le tableau suivant donne le crochet de Poisson entre les g\'en\'erateurs de $S^{A_2}$ :

$$\begin{array}{c|c|c|c|c|c|c|c}
&S_1&S_2&T_1&T_2&H&U_1&U_2\\
\hline
S_1&0&-H&0&U_1&-2S_1&2U_2&3T_1\\
\hline
S_2&H&0&-U_2&0&2S_2&-3T_2&-2U_1\\
\hline
T_1&0&U_2&0&-6S_1S_2+3H^2&-3T_1&-6S_1H&6S_1^2\\
\hline
T_2&-U_1&0&6S_1S_2-3H^2&0&3T_2&-6S_2^2&6S_2H\\
\hline
H&2S_1&-2S_2&3T_1&-3T_2&0&-U_1&U_2\\
\hline
U_1&-2U_2&3T_2&6S_1H&6S_2^2&U_1&0&-30S_1S_2+3H^2\\
\hline
U_2&-3T_1&2U_1&-6S_1^2&-6S_2H&-U_2&30S_1S_2-3H^2&0
\end{array}$$

En particulier, $H=\{S_2,S_1\}$, $U_1=\{S_2,T_2\}$ et $U_2=\{T_1,S_2\}$.\\

\parags Enfin, on peut ais\'ement montrer que les formules pr\'ec\'edentes
d\'efinissent un crochet de Poisson sur l'alg\`ebre de polyn\^omes 
$\C[S_1,,T_1,U_1,S_2,T_2,U_2,H]$. (Il suffit de v\'erifier l'identit\'e
de Jacobi sur les sept g\'en\'erateurs, ce qui se fait par un calcul direct
long, mais sans difficult\'e).

\subsection{Pr\'esentations de $S^{G_2}$ et $S^{B_2}$}

\parags Reprenons les notations du paragraphe \ref{sectG2}.
Comme $S^{G_2}$ est la somme directe des composantes homog\`enes de degr\'e pair
de $S^{A_2}$, on d\'eduit imm\'ediatement les r\'esultats suivants :
\begin{theo}
\begin{description}
\item[\it i)] $S^{G_2}$ est engendr\'ee par $S_1$, $T_1'=T_1^2$, $S_2$, $T_2'=T_2^2$, $Z=T_1T_2$,
$H$, $U_{11}=T_1U_1$, $U_{12}=T_1U_2$, $U_{21}=T_2U_1$ et $U_{22}=T_2U_2$.
\item[\it ii)] Posons $R'=\C[S_1,T'_1,S_2,T'_2]$. Tout \'el\'ement
de $S^{G_2}$ s'\'ecrit de mani\`ere unique :
\begin{eqnarray*}
Q&=&P_0+P_1H+P_2H^2+P_3H^3\\
&&+Q_0Z+Q_1ZH+Q_2ZH^2+Q_3ZH^3\\
&&+Q_{11}U_{11}+Q_{21}U_{21}+Q_{12}U_{12}+Q_{22}U_{22},
\end{eqnarray*}
o\`u les $P_i$, les $Q_j$ et les $Q_{kl}$ sont des \'el\'ements de $R'$.
\item[\it iii)] $S^{G_2}$ est engendr\'ee par 
$S_1$, $T_1'$, $S_2$, $T_2'$, $Z$,
$H$, $U_{11}$, $U_{12}$, $U_{21}$ et $U_{22}$ et les relations suivantes :
\begin{eqnarray*}
Z^2&=&T_1'T_2',\\
H^4&=&-4S_1^2S_2^2-3ZH+5S_1S_2H^2-S_2U_{11}-S_1U_{22},\\
HU_{11}&=&-3S_1Z-S_2U_{12},\\
HU_{12}&=&-3S_2T'_1-S_2U_{11},\\
HU_{21}&=&-3S_1T_2'-S_2U_{22},\\
HU_{22}&=&-3S_2Z-S_2U_{21},\\
U_{11}^2&=&12S_1S_2^2T'_1-3S_2T'_1H^2+3T'_1U_{22},\\
U_{11}U_{12}&=&9ZT'_1-12S_1S_2T'_1H+3T'_1H^3,\\
U_{11}U_{21}&=&12S_1S_2^2Z-3S_2ZH^2+3ZU_{22},\\
U_{11}U_{22}&=&9Z^2-12S_1S_2ZH+3ZH^3,\\
U_{12}^2&=&12S_1^2S_2T'_1-3S_1T'_1H^2+3T'_1U_{11},\\
U_{12}U_{21}&=&9Z^2-12S_1S_2ZH+3ZH^3,
\end{eqnarray*}
\begin{eqnarray*}
U_{12}U_{22}&=&12S_1^2S_2Z-3S_1ZH^2+3ZU_{11},\\
U_{21}^2&=&12S_1S_2^2T'_2-3S_2T'_2H^2+3T'_2U_{22},\\
U_{21}U_{22}&=&9ZT'_2-12S_1S_2T'_2H+3T'_2H^3,\\
U_{22}^2&=&12S_1^2S_2T'_2-3S_1T'_2H^2+3T'_2U_{11},\\
ZU_{11}&=&T'_1U_{21},\\
ZU_{12}&=&T'_1U_{22},\\
ZU_{21}&=&T'_2U_{11},\\
ZU_{22}&=&T'_2U_{12}.
\end{eqnarray*}
\end{description}
\end{theo}

\parags Reprenons les notations du paragraphe \ref{sectB2}. Par une m\'ethode semblable
\`a celle utilis\'ee pour $A_2$, on montrerait le r\'esultat suivant :

\begin{theo}
\begin{description}
\item[\it i)] $S^{B_2}$ est engendr\'ee par 
$S_1=x_1^2+x_2^2$, $S_2=x_1^2x_2^2$, 
$T_1=y_1^2+y_2^2$, $T_2=y_1^2y_2^2$, 
$Z_1=x_1y_1+x_2y_2$, $Z_2=x_1y_1x_2y_2$, $Z_3=x_1y_1^3+x_2y_2^3$
et $Z_4=x_1^3y_1+x_2^3y_3$.
\item[\it ii)] 
Posons $R=\C[S_1,T_1,S_2,T_2]$.
Il s'agit d'une alg\`ebre de polyn\^omes \`a quatre variables.
Tout \'el\'ement de $S^{B_2}$ s'\'ecrit de mani\`ere unique :
$$Q=P_0+P_1Z_1+P_2Z_1^2+P_3Z_1^3+P_4Z_1^4
+Q_2Z_2+Q_3Z_3+Q_4Z_4,$$
o\`u les $P_i$ et les $Q_j$ sont des \'el\'ements de $R$.
\item[\it iii)] $S^{B_2}$ est engendr\'ee par 
$S_1$, $T_1$, $S_2$, $T_2$, $Z_1$, $Z_2$, $Z_3$ et $Z_4$ 
et les relations suivantes :
\begin{eqnarray*}
Z_1^5&=&16S_2T_2-5S_2T_1^2-5S_1^2T_2+\frac{S_1^2T_1^2}{2}+\frac{3S_1T_1}{2}Z_1^2\\
&&+\left(4S_2T_1-\frac{S_1^2T_1}{2}\right)Z_3
+\left(4S_1T_2-\frac{S_1T_1^2}{2}\right)Z_4,\\
Z_1Z_2&=&\frac{S_1T_1}{4}Z_1+\frac{1}{4}Z_1^3-\frac{S_1}{4}Z_3-\frac{T_1}{4}Z_4,\\
Z_1Z_3&=&-S_1T_2+T_1Z_1^2-T_1Z_2,\\
Z_1Z_4&=&-S_2T_1+S_1Z_1^2-S_1Z_2,\\
Z_2^2&=&S_2T_2,\\
Z_2Z_3&=&\left(\frac{S_1T_1^2}{4}-S_1T_2\right)Z_1+\frac{T_1}{4}Z_1^3-\frac{S_1T_1}{4}Z_3
+\left(T_2-\frac{T_1^2}{4}\right)Z_4,\\
Z_2Z_4&=&\left(\frac{S_1^2T_1}{4}-S_2T_1\right)Z_1+\frac{S_1}{4}Z_1^3
+\left(S_2-\frac{S_1^2}{4}\right)Z_3-\frac{S_1T_1}{4}Z_4,\\
Z_3^2&=&-S_1T_1T_2+(T_1^2-T_2)Z_1^2+(4T_2-2T_1^2)Z_2,\\
Z_3Z_4&=&4S_2T_2-\frac{5S_2T_1^2}{4}-\frac{5S_1^2T_2}{4} 
+\frac{5S_1T_1}{4}Z_1^2-\frac{1}{4}Z_1^4-\frac{3S_1T_1}{2}Z_2,\\
Z_4^2&=&-S_1S_2T_1+(S_1^2-S_2)Z_1^2+(4S_2-2S_1^2)Z_2.
\end{eqnarray*}
\end{description}
\end{theo}

%% file: chap8.tex
\section{Sous-groupes de $(\pm 1)^n$}

\parag Soit maintenant $G$ un sous-groupe de $(\pm 1)^n$
agissant sur $S=\C[x_1,\ldots,x_n,y_1,\ldots,y_n]$ 
et sur l'alg\`ebre de Weyl $A=A_n(\C)$ de la mani\`ere suivante :
en notant $p_1,\ldots,p_n$, $q_1,\ldots,q_n$ les coordonn\'ees de $A_n(\C)$,
$$
\begin{array}{rcl}
(\epsilon_1,\ldots,\epsilon_n).x_i&=&\epsilon_i x_i,\\
(\epsilon_1,\ldots,\epsilon_n).y_i&=&\epsilon_i y_i\:;
\end{array}
\hspace{1cm}
\begin{array}{rcl}
(\epsilon_1,\ldots,\epsilon_n).p_i&=&\epsilon_i p_i,\\
(\epsilon_1,\ldots,\epsilon_n).q_i&=&\epsilon_i q_i.
\end{array}$$
\parag Nous avons le th\'eor\`eme suivant :
\begin{theo}
\label{theo12}
\begin{eqnarray*}
dim_\C \:HP_0(S^G) &=& \left\{\begin{array}{rl}
 0 &\mbox{s'il existe $i$ tel que $G \subseteq (\pm 1)^{i-1} \times \{1\}
\times (\pm1)^{n-i}$,} \\
  1&\mbox{sinon.}
\end{array}\right. \\
\\
dim_\C \:HH_0(A^G) &=& \left\{\begin{array}{rl}
 0 &\mbox{si $(-1,\ldots,-1) \notin G$,} \\
  1&\mbox{sinon.}
\end{array}\right. 
\end{eqnarray*}
\end{theo}

{\bf Preuve.} Pour tout $1\leq i\leq n$, $\g_i=Vect(x_i^2,x_iy_i,x_i^2)$ est inclus dans $S^G$
est forme une alg\`ebre de Lie isomorphe \`a $\sl$.
Par suite, comme dans les exemples pr\'ec\'edents, avec $\g=\g_1\oplus \ldots \oplus \g_n$,
la composante isotypique triviale $S^G_\g$ de $S^G$ sous l'action de $\g$ est
un compl\'ementaire de $\{S^G,S^G\}$. D'autre part :
\begin{eqnarray*}
S^G_\g&\subseteq& \bigcap_{i=1}^n Ker(\{x_i^2,.\}) \cap \bigcap_{i=1}^n Ker(\{y_i^2,.\})\\
&\subseteq&\bigcap_{i=1}^n Ker\left(x_i \frac{\partial}{\partial y_i}\right) 
\cap \bigcap_{i=1}^n Ker\left(y_i \frac{\partial}{\partial x_i}\right) \\
&\subseteq&\bigcap_{i=1}^n Ker\left(\frac{\partial}{\partial y_i}\right) 
\cap \bigcap_{i=1}^n Ker\left(\frac{\partial}{\partial x_i}\right) \\
&\subseteq & \C[x_1,\ldots,x_n] \cap \C[y_1,\ldots,y_n]\\
&\subseteq&\C.
\end{eqnarray*}
Donc $S^G_\g=\C$ et $S^G=\C + \{S^G,S^G\}$.
Par suite, $dim_\C \:HP_0(S^G)P=0$ si $1 \in \{S^G,S^G\}$ et $0$ sinon.

Supposons qu'il existe $i$ tel que $G \subseteq (\pm 1)^{i-1} \times \{1\}
\times (\pm1)^{n-i}$. Alors $x_i$ et $y_i \in S^G$ et donc $1=\{x_i,y_i\} \in \{S^G,S^G\}$.
Sinon, $G$ agit sans points fixes sur $Vect(x_1,\ldots,x_n,y_1,\ldots,y_n)$ et donc
la composante homog\`ene de $S^G$ de degr\'e $1$ est nulle. Comme $\{,\}$ est homog\`ene de degr\'e
$-2$, $1 \notin \{S^G,S^G\}$.

Enfin, rappelons que $dim_\C \:HH_0(A^G)$ est le nombre de classes de conjugaison 
de $G$ agissant sur $Vect(x_1,\ldots,x_n)$ sans points fixes. Le seul \'el\'ement
de $(\pm 1)^n$ agissant sans points fixes \'etant $(-1,\ldots,-1)$, le r\'esultat
sur $dim_\C \:HH_0(A^G)$ est imm\'ediat. $\Box$ \\

\parag {\bf Remarque.} Dans tous les cas, on a $dim_\C \:HP_0(S^G)\geq dim_\C \:HH_0(A^G)$.
Ces deux nombres ne sont pas n\'ecessairement \'egaux, comme le montre le troisi\`eme exemple 
ci-dessous :
\begin{description}
\item[\it i)] Pour $G=(\pm1)^{n-1}\times (1)$, alors $dim_\C \:HP_0(S^G)=dim_\C \:HH_0(A^G)=0$.
\item[\it ii)] Pour $G=\{(\epsilon_1,\ldots,\epsilon_n)\in (\pm1)^n\:/\:
\epsilon_1\ldots\epsilon_n=1\}$, avec $n$ pair, $dim_\C \:HP_0(S^G)=dim_\C \:HH_0(A^G)=1$.
\item[\it iii)] Pour $G=\{(\epsilon_1,\ldots,\epsilon_n)\in (\pm1)^n\:/\:
\epsilon_1\ldots\epsilon_n=1\}$, avec $n$ impair sup\'erieur ou \'egal \`a $3$, 
$dim_\C \:HP_0(S^G)=1$ et $dim_\C \:HH_0(A^G)=0$.\\
\end{description}

\parag En particulier, pour l'exemple $iii)$, les constantes sont dans $[A^G,A^G]$. Voici
une mani\`ere de les exprimer :
\begin{prop}
Dans $A^G$, avec les notations de l'exemple {\it iii)} ci-dessus :
$$\sum_{(x_i,y_i)=(p_i,q_i)\mbox{ ou }(q_i,p_i)} 
(-1)^{\alpha_{x,y}} [x_1\ldots x_n,y_1\ldots y_n]=2,$$
o\`u $\alpha_{x,y}$ d\'esigne le nombre de $x_i$ \'egaux \`a $q_i$.
\end{prop}
Par exemple, pour $n=1$ et pour $n=3$ :
\begin{eqnarray*}
2&=&[p_1,q_1]-[q_1,p_1],\\
2&=&[p_1p_2p_3,q_1q_2q_3]-[q_1p_2p_3,p_1q_2q_3]-[p_1q_2p_3,q_1p_2q_3]-[p_1p_2q_3,q_1q_2p_3]\\
&&+[p_1q_2q_3,q_1p_2p_3]+[q_1p_2q_3,p_1q_2p_3]+[q_1q_2p_3,p_1p_2q_3]-[q_1q_2q_3,p_1p_2p_3].
\end{eqnarray*}

{\bf Preuve.} 

{\it Premi\`ere \'etape.} Dans $A$, on pose, pour tout $1\leq i\leq n$ :
$$E_i=\frac{p_i^2}{2},\: F_i=-\frac{q_i^2}{2},\: H_i=-\frac{p_iq_i+q_ip_i}{2}.$$
Alors $\g_i=Vect(E_i,F_i,H_i)$ munie du crochet $[,]$ est une sous-alg\`ebre Lie de $A$
isomorphe \`a $\sl$. En posant $\g=\g_1+\ldots+\g_n$, $\g$ est une sous-alg\`ebre de Lie
de $A$ isomorphe \`a $\sl^{\oplus n}$. Elle est de plus incluse dans $A^G$, donc $A^G$
est ainsi muni d'une structure de $\g$-module. 
De mani\`ere semblable \`a la preuve du th\'eor\`eme \ref{theo12}, la composante isotypique 
triviale de $A^G$ est $\C$. 

{\it Deuxi\`eme \'etape.} Consid\'erons l'\'el\'ement suivant de l'alg\`ebre $A \otimes A$ :
$$X=(p_1\otimes q_1-q_1 \otimes p_1) \ldots (p_n\otimes q_n-q_n \otimes p_n).$$
En d\'eveloppant :
$$X=\sum_{(x_i,y_i)=(p_i,q_i)\mbox{ ou }(q_i,p_i)}
 (-1)^{\alpha_{x,y}} x_1\ldots x_n \otimes y_1\ldots y_n,$$
donc $X \in A^G \otimes A^G$. 
D'autre part, $\g$ agit par d\'erivation sur $A$ 
et donc agit \'egalement par d\'erivation sur $A\otimes A$.
Par suite, pour tout $1\leq i \leq n$ :
\begin{eqnarray*}
E_i.X&=&(p_1\otimes q_1-q_1 \otimes p_1) \ldots E_i.(p_i\otimes q_i-q_i \otimes p_i)
\ldots (p_n\otimes q_n-q_n \otimes p_n)+0\\
&=&(p_1\otimes q_1-q_1 \otimes p_1) \ldots (E_i.p_i\otimes q_i+p_i\otimes E_i.q_i
-E_i.q_i \otimes p_i-q_i \otimes E_i.p_i)
\ldots (p_n\otimes q_n-q_n \otimes p_n)\\
&=&(p_1\otimes q_1-q_1 \otimes p_1) \ldots (0+p_i\otimes p_i
-p_i \otimes p_i-0)
\ldots (p_n\otimes q_n-q_n \otimes p_n)\\
&=&0,\\
\\
F_i.X&=&(p_1\otimes q_1-q_1 \otimes p_1) \ldots F_i.(p_i\otimes q_i-q_i \otimes p_i)
\ldots (p_n\otimes q_n-q_n \otimes p_n)+0\\
&=&(p_1\otimes q_1-q_1 \otimes p_1) \ldots (F_i.p_i\otimes q_i+p_i\otimes F_i.q_i
-F_i.q_i \otimes p_i-q_i \otimes F_i.p_i)
\ldots (p_n\otimes q_n-q_n \otimes p_n)\\
&=&(p_1\otimes q_1-q_1 \otimes p_1) \ldots (q_i\otimes q_i+0
-0-q_i \otimes q_i)
\ldots (p_n\otimes q_n-q_n \otimes p_n)\\
&=&0,\\
\\
H_i.X&=&(p_1\otimes q_1-q_1 \otimes p_1) \ldots H_i.(p_i\otimes q_i-q_i \otimes p_i)
\ldots (p_n\otimes q_n-q_n \otimes p_n)+0\\
&=&(p_1\otimes q_1-q_1 \otimes p_1) \ldots (H_i.p_i\otimes q_i+p_i\otimes H_i.q_i
-H_i.q_i \otimes p_i-q_i \otimes H_i.p_i)
\ldots (p_n\otimes q_n-q_n \otimes p_n)\\
&=&(p_1\otimes q_1-q_1 \otimes p_1) \ldots (p_i\otimes q_i-p_i\otimes q_i
+q_i \otimes p_i-q_i \otimes p_i)
\ldots (p_n\otimes q_n-q_n \otimes p_n)\\
&=&0,
\end{eqnarray*}
donc $X$ est dans la composante isotypique triviale de $A^G \otimes A^G$.\\

{\it Troisi\`eme \'etape.} L'identit\'e de Jacobi montre que 
$[,]:A^G \otimes A^G \longrightarrow A^G$
est un morphisme de $\g$-modules. Par suite, $[X]$ est dans la composante isotypique 
triviale de $A^G$, c'est-\`a-dire dans $\C$. Pour calculer $[X]$, il suffit donc de
calculer le terme constant de chacun des $[x_1\ldots x_n,y_1\ldots y_n]$ dans la base
$(p_1^{\alpha_1}\ldots p_n^{\alpha_n}q_1^{\beta_1}\ldots q_n^{\beta_n})_
{\alpha_i,\beta_i \in \mathbb{N}}$ de $A$. Trois cas se pr\'esentent.
\begin{description}
\item[\it i)] Il existe $i,j\in \{1,\ldots,n\}$ tel que $x_i=p_i$ et $y_j=p_j$. Alors :
$$[x_1\ldots x_n,y_1\ldots y_n]=\underbrace{x_1\ldots x_n y_1\ldots y_n}_{\in p_iA}
-\underbrace{y_1\ldots y_nx_1\ldots x_n}_{\in p_jA}\in p_iA+p_jA.$$
Or une base de $p_iA+p_jA$ est 
$(p_1^{\alpha_1}\ldots p_n^{\alpha_n}q_1^{\beta_1}\ldots q_n^{\beta_n})_
{\alpha_i>0 \mbox{ ou }\alpha_j>0}$. Par suite le terme constant de $[x_1\ldots x_n,y_1\ldots y_n]$
est nul.
\item[\it ii)] Tous les $x_i$ sont \'egaux \`a $p_i$. Une r\'ecurrence montre que pour
tout $k \in \mathbb{N}^*$, dans l'alg\`ebre de Weyl $A_k(\C)$ :
$$[p_1\ldots p_k,q_1\ldots q_k]=\sum_{j=0}^{k-1} \sum_{1\leq i_1<\ldots<i_j\leq k}
(-1)^{k-j-1}p_{i_1}\ldots p_{i_j} q_{i_1}\ldots q_{i_j}.$$
Par suite, le terme constant de $[p_1\ldots p_n,q_1\ldots q_n]$ est $(-1)^{n-1-0}=1$
car $n$ est impair.
\item[\it iii)] Tous les $x_i$ sont \'egaux \`a $q_i$. Par antisym\'etrie de $[,]$,
le terme constant de $[q_1\ldots q_n,p_1\ldots p_n]$ est $-1$. 
\end{description}
En conclusion, le terme constant de $[X]$ est $1+(-1)^n(-1)=2$. Comme $[X]\in\C$,
$[X]=2$. $\Box$\\

\parag {\bf Remarque.} On a montr\'e au passage que $S$ et $A$ sont deux $\sl^{\oplus n}$-modules.
On peut montrer que le morphisme de sym\'etrisation
de $S$ dans $A$ est un isomorphisme de $\sl^{\oplus n}$-modules ;
comme il commute \`a l'action de $G$, il induit un isomorphisme de 
$\sl^{\oplus n}$-modules de $S^G$ dans $A^G$.

%% file: chap9.tex
\section{Sous-groupe de $\left(\frac{\mathbb{Z}}{3\mathbb{Z}}\right)^n$}

\parag 
Soit $n \in \mathbb{N}$, plus grand que $2$.
Consid\'erons le sous-groupe suivant de $\left(\frac{\mathbb{Z}}{3\mathbb{Z}}\right)^n$ :
$$G=\left\{(\overline{k}_1,\ldots,\overline{k}_n)\in \left(\frac{\mathbb{Z}}{3\mathbb{Z}}\right)^n
\:/\: \overline{k}_1+\ldots+\overline{k}_n=\overline{0}\right\}.$$
Ce groupe agit sur $S=\C[x_1,\ldots,x_n,y_1,\ldots,y_n]$ munie du crochet de Poisson standard
et sur l'alg\`ebre de Weyl $A=A_n(\C)$ de la mani\`ere suivante :
en notant $p_1,\ldots,p_n$, $q_1,\ldots,q_n$ les coordonn\'ees de $A_n(\C)$,
\begin{eqnarray*}
(\overline{k}_1,\ldots,\overline{k}_n).x_i&=&\zeta^{k_i} x_i,\\
(\overline{k}_1,\ldots,\overline{k}_n).y_i&=&\zeta^{-k_i} y_i\:;\\
(\overline{k}_1,\ldots,\overline{k}_n).p_i&=&\zeta^{k_i} p_i,\\
(\overline{k}_1,\ldots,\overline{k}_n).q_i&=&\zeta^{-k_i} q_i,
\end{eqnarray*}
o\`u $\zeta=e^{\frac{2i\pi}{3}}$.\\

\parag Nous avons le th\'eor\`eme suivant :
\begin{theo}
\begin{eqnarray*}
dim_\C \:HP_0(S^G) &=&2^n-2,\\
dim_\C \:HH_0(A^G) &=&\frac{2}{3}\left(2^{n-1}-(-1)^{n-1}\right). 
\end{eqnarray*}
\end{theo}

{\bf Preuve.} 

{\it Premi\`ere \'etape.}
Pour tout $i \in \{1,\ldots,n\}$, 
consid\'erons $t_i=x_iy_i \in S^G$ et posons $\g=Vect(t_1,\ldots,t_n)$ ;
il s'agit d'une sous-alg\`ebre de Lie de $S^G$. De plus :
$$\{t_i,x_1^{\alpha_1}\ldots x_n^{\alpha_n}
y_1^{\beta_1}\ldots y_n^{\beta_n}\}=(\beta_i-\alpha_i)x_1^{\alpha_1}\ldots x_n^{\alpha_n}
y_1^{\beta_1}\ldots y_n^{\beta_n}.$$
Donc $\g$ agit de mani\`ere semi-simple sur $S$ et donc aussi sur $S^G$.
Sa composante isotypique triviale est $\C[t_1,\ldots,t_n]$.
Comme dans les exemples pr\'ec\'edents, on en d\'eduit :
\begin{eqnarray*}
S^G&=&\C[t_1,\ldots,t_n]+\{S^G,S^G\},\\
HP_0(S^G)&=&\frac{\C[t_1,\ldots,t_n]}{\C[t_1,\ldots,t_n]\cap\{S^G,S^G\}}.
\end{eqnarray*}
De plus, dans $S^G$, si $(\alpha_1,\ldots,\alpha_n) \in \mathbb{N}^n$ :
\begin{eqnarray*}
\{t_1^{\alpha_1}\ldots t_n^{\alpha_n}x_i^3,y_i^3\}
&=&\{x_1^{\alpha_1} \ldots x_i^{\alpha_i+3} \ldots x_n^{\alpha_n},y_i^3\}
y_1^{\alpha_1}\ldots y_n^{\alpha_n}\\
&=&3(\alpha_i+3)x_1^{\alpha_1} \ldots x_i^{\alpha_i+2} \ldots x_n^{\alpha_n}
y_1^{\alpha_1}\ldots y_i^{\alpha_i+2} \ldots y_n^{\alpha_n}\\
&=&t_1^{\alpha_1}\ldots t_i^{\alpha_i+2}\ldots t_n^{\alpha_n}.
\end{eqnarray*}
Donc $<t_1^2,\ldots,t_n^2> \subseteq \C[t_1,\ldots,t_n]\cap\{S^G,S^G\}$. Par suite :
$$HP_0(S^G)=\frac{\displaystyle \frac{\C[t_1,\ldots,t_n]}{<t_1^2,\ldots,t_n^2>}}
{\displaystyle\frac{\C[t_1,\ldots,t_n]\cap\{S^G,S^G\}}{<t_1^2,\ldots,t_n^2>}}.$$
On pose $R=\frac{\displaystyle \C[t_1,\ldots,t_n]}{\displaystyle<t_1^2,\ldots,t_n^2>}$.
Il s'agit d'une alg\`ebre ayant pour base $(\overline{t_1^{\alpha_1}\ldots t_n^{\alpha_n}})
_{0\leq \alpha_i \leq 2}$, donc de dimension $2^n$.
On pose $J=\displaystyle\frac{\C[t_1,\ldots,t_n]\cap\{S^G,S^G\}}{<t_1^2,\ldots,t_n^2>}$,
de sorte que $HP_0(S^G)=R/J$.\\

{\it Deuxi\`eme \'etape.} $S^G$ est engendr\'ee par les mon\^omes fix\'es par $G$. 
Soient $(\alpha_1,\ldots,\alpha_n)$ 
et $(\beta_1,\ldots,\beta_n) \in \mathbb{N}^n$.
Alors :\\

\begin{tabular}{cl}
&$x_1^{\alpha_1}\ldots x_n^{\alpha_n}
y_1^{\beta_1}\ldots y_n^{\beta_n} \in S^G$\\[2mm]
$\Longleftrightarrow$&
$\forall (\overline{k}_1,\ldots,\overline{k}_n)\in G$,
$k_1(\alpha_1-\beta_1)+\ldots +k_n(\alpha_n-\beta_n)\equiv 0[3]$\\[2mm]
$\Longleftrightarrow$&
$\exists \lambda \in \{0,1,2\}$, $\forall i \in \{1,\ldots,n\}$, 
$\alpha_i-\beta_i \equiv \lambda [3]$.
\end{tabular}\\

On a donc :
\begin{eqnarray*}
&&\C[t_1,\ldots,t_n]\cap\{S^G,S^G\}\\
&=&Vect\left(\begin{array}{c}
\{x_1^{\alpha_1}\ldots x_n^{\alpha_n} y_1^{\beta_1}\ldots y_n^{\beta_n},
 x_1^{\gamma_1}\ldots x_n^{\gamma_n} y_1^{\delta_1}\ldots y_n^{\delta_n}\}\:/\:\\
  x_1^{\alpha_1}\ldots x_n^{\alpha_n} y_1^{\beta_1}\ldots y_n^{\beta_n},
 x_1^{\gamma_1}\ldots x_n^{\gamma_n} y_1^{\delta_1}\ldots y_n^{\delta_n} \in S^G,\:
\forall i,\: \alpha_i-\beta_i+\gamma_i-\delta_i=0
\end{array}\right)\\
&=&Vect\left(\begin{array}{c}
\displaystyle \sum_{i=1}^n (\alpha_i \gamma_i-\beta_i\delta_i)t_1^{\alpha_1+\gamma_1}
\ldots t_i^{\alpha_i+\gamma_i-1}\ldots t_n^{\alpha_n+\gamma_n}\:/\:\\
\exists 0\leq \lambda\leq 2,\: \forall i, \: \alpha_i-\beta_i\equiv\lambda[3],\:
\alpha_i-\beta_i+\gamma_i-\delta_i=0
\end{array}\right).
\end{eqnarray*}
Supposons qu'il existe $i$, tel que $\alpha_i$, $\beta_i$, $\gamma_i$ 
ou $\delta_i$ soit sup\'erieur ou \'egal \`a $3$. Alors tous les mon\^omes apparaissant
dans $\{x_1^{\alpha_1}\ldots x_n^{\alpha_n} y_1^{\beta_1}\ldots y_n^{\beta_n},
 x_1^{\gamma_1}\ldots x_n^{\gamma_n} y_1^{\delta_1}\ldots y_n^{\delta_n}\}$
 sont dans $<t_1^2,\ldots,t_n^2>$. Par suite :
$$J=Vect\left(\begin{array}{c}
\displaystyle \sum_{i=1}^n (\alpha_i \gamma_i-\beta_i\delta_i)\overline{t_1^{\alpha_1+\gamma_1}
\ldots t_i^{\alpha_i+\gamma_i-1}\ldots t_n^{\alpha_n+\gamma_n}}\:/\:\\
\exists 0\leq \lambda\leq 2,\: \forall i, \: \alpha_i-\beta_i\equiv\lambda[3],\:
\alpha_i-\beta_i+\gamma_i-\delta_i=0,\: 0\leq \alpha_i,\beta_i,\gamma_i,\delta_i\leq 2
\end{array}\right).$$
Par la suite, nous noterons $J$ cet espace.
Consid\'erons un \'el\'ement g\'en\'erateur de $J$
pour lequel $\lambda=0$. Comme pour tout $i$,
$\alpha_i-\beta_i\equiv 0[3]$ et $-2\leq \alpha_i-\beta_i \leq 2$, 
n\'ecessairement $\alpha_i=\beta_i$. De m\^eme, $\gamma_i=\delta_i$, donc
$(\alpha_i \gamma_i-\beta_i\delta_i)=0$. Par suite :
$$J=Vect\left(\begin{array}{c}
\displaystyle \sum_{i=1}^n (\alpha_i \gamma_i-\beta_i\delta_i)\overline{t_1^{\alpha_1+\gamma_1}
\ldots t_i^{\alpha_i+\gamma_i-1}\ldots t_n^{\alpha_n+\gamma_n}}\:/\:\\
\exists 1\leq \lambda\leq 2,\: \forall i, \: \alpha_i-\beta_i\equiv\lambda[3],\:
\alpha_i-\beta_i+\gamma_i-\delta_i=0,\: 0\leq \alpha_i,\beta_i,\gamma_i,\delta_i\leq 2
\end{array}\right).$$
\\

{\it Troisi\`eme \'etape.}
Consid\'erons un \'el\'ement g\'en\'erateur de $J$
pour lequel $\lambda=1$. Comme pour tout $i$,
$\alpha_i-\beta_i\equiv 1[3]$ et $-2\leq \alpha_i-\beta_i \leq 2$, 
n\'ecessairement :
$$\left\{\begin{array}{c}\alpha_i=\beta_i+1,\\
\gamma_i=\delta_i-1
\end{array}\right.
\mbox{ ou }
\left\{\begin{array}{c}\alpha_i=\beta_i-2,\\
\gamma_i=\delta_i+2
\end{array}\right.
$$ 
\begin{description}
\item[\it a)] Dans le premier cas, on obtient $(\alpha_i,\beta_i, \gamma_i,\delta_i)
\in \{(1,0,0,1),(1,0,1,2),(2,1,0,1),(2,1,1,2)\}$. Quatre sous-cas apparaissent :
\begin{description}
\item[\it i)] Il existe $j$ tel que $(\alpha_j,\beta_j, \gamma_j,\delta_j)= (2,1,1,2)$.
Alors tous les mon\^omes apparaissant dans le g\'en\'erateur ont une puissance en $t_j$ 
sup\'erieur ou \'egale \`a $2$, donc sont dans $<t_1^2,\ldots,t_n^2>$ :  le g\'en\'erateur est nul.
\item[\it ii)] Il existe $j,k$ distincts tels que 
$(\alpha_j,\beta_j, \gamma_j,\delta_j)$  et $(\alpha_k,\beta_k, \gamma_k,\delta_k)$ 
sont des \'el\'ements de $\{(1,0,1,2),(2,1,0,1)\}$. Alors tous les mon\^omes
apparaissant dans le g\'en\'erateur ont une puissance en $t_j$ ou en $t_k$ 
sup\'erieur ou \'egale \`a  $2$, donc sont dans $<t_1^2,\ldots,t_n^2>$ :  le g\'en\'erateur est nul.
\item[\it iii)] Il existe $j$ tel que 
$(\alpha_j,\beta_j, \gamma_j,\delta_j)\in \{(1,0,1,2),(2,1,0,1)\}$ et si
$k\neq j$, $(\alpha_k,\beta_k, \gamma_k,\delta_k)=(1,0,0,1)$. 
Alors le g\'en\'erateur vaut $2 \overline{t_1\ldots t_n}$.
\item[\it iv)] Si pour tout $j$, 
$(\alpha_j,\beta_j, \gamma_j,\delta_j)= (1,0,0,1)$, alors le g\'en\'erateur vaut
$\displaystyle \sum_{i=1}^n \overline{t_1\ldots t_{i-1}t_{i+1}\ldots t_n}$.
\end{description}
\item[\it b)] Dans le second cas, pour tout $i$, 
$(\alpha_i,\beta_i, \gamma_i,\delta_i)=(0,2,2,0)$ et donc tous les mon\^omes
apparaissant dans le g\'en\'erateur ont une puissance en $t_1$ ou en $t_2$ 
sup\'erieur ou \'egal \`a  $2$, donc sont dans $<t_1^2,\ldots,t_n^2>$ : le g\'en\'erateur est nul.
\end{description}
Le cas o\`u $\lambda=2$ se traite de la m\^eme mani\`ere, en permutant
les r\^oles de $(\alpha_i,\beta_i)$ et $(\gamma_i,\delta_i)$. On obtient donc :
$$J=Vect\left(\overline{t_1\ldots t_n},
\sum_{i=1}^n \overline{t_1\ldots t_{i-1}t_{i+1}\ldots t_n}
\right).$$
Il s'agit donc d'un espace de dimension $2$.
Comme $HP_0(S^G)=R/J$, $HP_0(S^G)$ est de dimension $2^n-2$.\\

{\it Quatri\`eme \'etape.} La dimension de $HH_0(A^G)$ est le nombre de classes de conjugaison
de $G$ agissant sans points fixes, c'est-\`a-dire le nombre d'\'el\'ements de $G$
sans composantes nulles. On pose, pour $n \geq 1$ :
$$a_n=card\left(\left\{(\overline{k}_1,\ldots,\overline{k}_n) \in \left(
\frac{\mathbb{Z}}{3\mathbb{Z}}-\{\overline{0}\}\right)^n\:/\:
\overline{k}_1+\ldots+\overline{k}_n=\overline{0}\right\}\right).$$
On pose \'egalement $a_0=1$.
On a alors :
$$|G|=3^{n-1}=1+\sum_{i=1}\binom{n}{i}a_i.$$
Par suite, en passant aux s\'eries g\'en\'eratrices exponentielles :
$$\sum_{n=1}^{+\infty} \frac{3^{n-1}}{n!}x^n=
\left(\sum_{n=0}^{+\infty} \frac{a_n}{n!}x^n\right)
\left(\sum_{n=0}^{+\infty} \frac{x^n}{n!}\right)-1.$$
Par suite :
$$\frac{1}{3}e^{3x}-\frac{1}{3}=e^x\left(\sum_{n=0}^{+\infty} \frac{a_n}{n!}x^n\right)-1,$$
$$\frac{1}{3}e^{2x}+\frac{2e^{-x}}{3}=\sum_{n=0}^{+\infty} \frac{a_n}{n!}x^n.$$
Donc $\displaystyle a_n=dim_\C \: HH_0(A^G)=\frac{1}{3}\left(2^n+2(-1)^n \right)
=\frac{2}{3}\left(2^{n-1}-(-1)^{n-1} \right)$.
 $\Box$ \\

\parag {\bf Remarque.} Par suite, 
$$dim_\C \:HP_0(S^G)-dim_\C \:HH_0(A^G)=\frac{2}{3}
\left(2^n+(-1)^{n-1}-3 \right).$$
Cette diff\'erence tend donc vers $+\infty$ quand $n$ tend vers $+\infty$.

%% file: Poisson2.bbl
\providecommand{\bysame}{\leavevmode\hbox to3em{\hrulefill}\thinspace}
\providecommand{\MR}{\relax\ifhmode\unskip\space\fi MR }
% \MRhref is called by the amsart/book/proc definition of \MR.
\providecommand{\MRhref}[2]{%
  \href{http://www.ams.org/mathscinet-getitem?mr=#1}{#2}
}
\providecommand{\href}[2]{#2}
\begin{thebibliography}{AFLS00}

\bibitem[AF00]{AF}
Jacques \bgroup Alev\egroup{} et Daniel~R. \bgroup Farkas\egroup{},
  \emph{{F}inite group actions on {P}oisson algebras}, The orbit method in
  geometry and physics (Marseille, 2000), Progr. Math., vol. 213,
  Birkh{\"a}user Boston, 2000, pp.~9--28.

\bibitem[AFLS00]{AFLS}
Jacques \bgroup Alev\egroup{}, Marco~A. \bgroup Farinati\egroup{}, Thierry
  \bgroup Lambre\egroup{} et Andrea~L. \bgroup Solotar\egroup{},
  \emph{Homologie des invariants d'une alg\`ebre de {W}eyl sous l'action d'un
  groupe fini}, J. Algebra \textbf{232} (2000), no.~2, 564--577.

\bibitem[AL98]{AL}
Jacques \bgroup Alev\egroup{} et Thierry \bgroup Lambre\egroup{},
  \emph{Comparaison de l'homologie de {H}ochschild et de l'homologie de
  {P}oisson pour une d\'eformation des surfaces de {K}lein}, Algebra and
  operator theory (Tashkent, 1997), Kluwer Acad. Publ., Dordrecht, 1998,
  pp.~25--38.

\bibitem[BEG04]{BEG}
Yuri \bgroup Berest\egroup{}, Pavel \bgroup Etingof\egroup{} et Victor \bgroup
  Ginzburg\egroup{}, \emph{Morita equivalence of {C}herednik algebras}, J.
  Reine Angew. Math. \textbf{568} (2004), 81--98.

\bibitem[BG03]{BG}
Kenneth~A. \bgroup Brown\egroup{} et Iain \bgroup Gordon\egroup{},
  \emph{Poisson orders, symplectic reflection algebras and representation
  theory}, J. Reine Angew. Math. \textbf{559} (2003), 193--216.

\bibitem[EG02]{EG}
Pavel \bgroup Etingof\egroup{} et Victor \bgroup Ginzburg\egroup{},
  \emph{Symplectic reflection algebras, {C}alogero-{M}oser space, and deformed
  {H}arish-{C}handra homomorphism}, Invent. Math. \textbf{147} (2002), no.~2,
  243--348.

\bibitem[Fu05]{FU}
Baohua \bgroup Fu\egroup{}, \emph{A survey on symplectic singularities and
  resolutions}, math.AG/05\:10346, 2005.

\end{thebibliography}
